\documentclass[12pt]{amsart}
\usepackage{amssymb}
\newfont{\sheaf}{eusm10 scaled\magstep1}

\def\CC{{\mathbb{C}}}

\newcommand{\QED}{\hspace*{\fill}$Q.E.D.$}

\newcommand{\NNN}{\ensuremath{\mathcal{N}}}

\newcommand{\V}{\ensuremath{\mathcal{V}}}

\newcommand{\ra}{\ensuremath{\rightarrow}}
\newcommand{\SSS}{\ensuremath{\mathcal{S}}}
\newcommand{\CCC}{\ensuremath{\mathcal{C}}}
\def\eea{\end{eqnarray*}}
\def\bea{\begin{eqnarray*}}

\def\A{{\mathcal{A}}}
\def\B{{\mathcal{B}}}
\def\C{{\mathbb{C}}}
\def\E{{\mathcal{E}}}

\def\K{{\mathcal{K}}}
\def\L{{\mathcal{L}}}
\def\M{{\mathcal{M}}}

\def\P{{\mathcal{P}}}
\def\R{{\mathcal{R}}}
\def\T{{\mathcal{T}}}

\def\hol{{\mathcal{O}}}

\def\PP{{\mathbb{P}}}
\def\Proj{{\mathbf{Proj}}}
\def\qed{{\hfill Q.E.D.}}
\def\text{{\hbox}}
\newcommand{\Proof}{{\it Proof. }}
\newtheorem{teo}{Theorem}[section]
\newtheorem{df}[teo]{Definition}
\newtheorem{lem}[teo]{Lemma}

\newtheorem{oss}[teo]{Remark}
\newtheorem{prop}[teo]{Proposition}
\newtheorem{ex}[teo]{Example}

\begin{document}
\thanks{ The
research of the  authors was performed through the years in the realm  of the
   DFG SCHWERPUNKT "Globale Methode in der komplexen Geometrie",  of the
EAGER EEC Project, and of the VIGONI-DAAD exchange Program.
The second author is a member of G.N.S.A.G.A. of
I.N.d.A.M.  and gratefully acknowledges the hospitality
of the Institute of Mathematics Simion Stoilow of the Rumanian Academy of
Sciences, in march 2002.\\
A.M.S. Subject classifications: 14D06, 14J29, 11G30.}
\title{Fibrations of low genus, I.\\
Fibrations en courbes de genre petit, I.}
\author{Fabrizio Catanese -- Roberto Pignatelli}
\date{5 april 2005} \maketitle
\begin{abstract}
In the present paper we consider fibrations
$f: S \ra B$ of an algebraic surface
over a curve $B$, with general fibre a curve of genus $g$.
Our main results are:

1) A structure theorem for such fibrations in the case
where $g=2$

2) A structure theorem for such fibrations in the case
where $g=3$,
the general fibre is nonhyperelliptic, and each fibre is 2-connected

3) A theorem giving a complete description of the
moduli space  of minimal surfaces of general type with $  p_g = q=
1$, $K^2_S = 3,$
showing in particular that it has four unirational connected
components

4) Other applications of the two structure theorems.

\noindent R\'ESUM\'E.
Dans cet article nous consid\'erons des fibrations $f: S \ra B$ d'une surface
alg\'ebrique $S$ sur une courbe $B$, dont la fibre g\'en\'erale est 
une courbe de
genre $g$. Nos r\'esultats principaux sont les suivants :

1) Un th\'eor\`eme de structure pour de telles fibrations dans le cas
$g=2$

2) Un th\'eor\`eme de structure pour de telles fibrations dans le cas
ou $g=3$, la
fibre g\'en\'erale est non hyperelliptique, et chaque fibre est
2-connexe

3) Un th\'eor\`eme donnant une description compl\`ete de l'espace de modules
des surfaces minimales de type g\'en\'eral  avec $  p_g = q=
1$, $K^2_S = 3,$ en
montrant en particulier qu'il comporte quatre  composantes
connexes qui sont unirationnelles

4) D' autres applications de ces deux th\'eor\`emes de structure.
\end{abstract}
\tableofcontents
\section{Introduction}
The study of fibrations $ f : S \ra B$ of an algebraic surface $S$ over
a curve $B$ lies at the heart of surface classification (cf. e.g.
\cite{enr}, \cite{bpv}).

We denote  by $g$ the genus of a general fibre;
  any surface birational to $S$ is called a {\bf a ruled surface} if
  $g=0$, and {\bf an elliptic surface} if  $g=1$.
  Establishing the existence of such fibrations with $g=0,1$
plays a prominent role  in the Enriques classification of
algebraic surfaces.

Genus 2 fibrations play a special role for surfaces of general type: the
presence of a genus~2 fibration constitutes the so-called {\bf standard
exception} to the birationality of the bicanonical map (cf.\ \cite{bom},
\cite{ccml}, \cite{cil}).

 From a birational point of view, the datum of the fibration
$ f : S \ra B$ is equivalent to the datum of a curve $C$ of genus $g$
over the function field $\C (B)$ of $B$.

Apart from complications coming from the fact that $\C (B)$ is not
algebraically closed, curves of low genus $ g \leq 5$ are easily described
(mostly as
complete intersection curves in a projective space), so one can hope
to construct, resp. describe such fibrations as complete intersections
in some projective bundle over the base curve $B$.

To do this, one needs a fixed biregular model for the birational
equivalence class of $ f : S \ra B$.
The classical approach is to consider the so called
{\bf relative minimal model},
which is unique except in the case $g=0$ (cf. \cite{shaf}), where
a relative minimal model is a $\PP^1$ bundle, but  the so called
elementary transformations yield nonisomorphic surfaces out of the
same birational genus $0$ fibration.

When $g=0$, all the fibres of a relative minimal fibration are
iso\-morphic to $\PP^1$; for $g=1$  Kodaira gave, as a preliminary 
tool for his deep
investigations of elliptic surfaces, a short complete list of the possible
fibres (cf. \cite{kod}).
In the $g=2$ case  a similar but overlong list was provided by
Ogg in \cite{ogg}. For example, Bombieri
(\cite{bomman}) was able to use Ogg's classification to prove that the
genus of the bicanonical pencil of a numerical Godeaux surface cannot
be~2. In this vein, we demonstrate in this paper the power of our new
methods by providing a half-page proof of Bombieri's result.

To explain our new results, we first remark that
the existing literature on fibrations with low $g$ usually divides 
into two lines
of research:

(1) Papers devoted to the necessary numerical restrictions that must
  be satisfied by surfaces admitting such a fibration (see e.g. \cite{hor2},
\cite{xiao2}, \cite{kon1})

(2) Papers devoted to proving existence or devoted to the classification of
such fibrations.

For instance, the classification of surfaces with irregularity
$q=1$ leads naturally to the above issue (2), because the
Albanese map of such surfaces is a genus $g$ fibration over a curve
$B$ of genus $1$.
Even when restricting to the case $p_g= q=1$, where by general results
such surfaces belong to  a finite number of families, the classification
has turned out  to be quite a hard task.
The first result is quite old, and due independently to Bombieri and
the first named author and to  Horikawa.

\

{\bf  Theorem (\cite{k^2=2}, \cite{hor2})} {\it A minimal surface with $K^2_S
= 2, p_g= q=1$  is the double cover of the symmetric product $B^{(2)}$ of
the elliptic curve $ B : = Alb(S)$,
with a branch divisor $\Delta$ which belongs to a fixed linear system
of   6-sections of
the
$\PP^1$-bundle
$B^{(2)} \ra B$.

Their moduli space is unirational of dimension $7$.}

\

Recall in fact that
the symmetric product $B^{(h)} $ of  an elliptic curve parametrizes
effective divisors of degree $h$ on $B$, and, using the group structure
on $B$, we get the Abel-Jacobi map $ \pi : B^{(h)} \ra B$,
associating to the divisor $ (P_1, \dots ,P_h)$ the point
$ \pi (P_1, \dots ,P_h) : = P_1 + \dots  +P_h.$
Abel's theorem shows that $ \pi$ makes $B^{(h)} $ a $\PP^{h-1}$ bundle
over $B$, and in fact we get the projectivization
of an indecomposable vector bundle, denoted
$ E(h,1)$ by Atiyah (\cite{atiyah}).

Ciliberto later proposed to the first author to construct surfaces with
$p_g= q=1$ and higher $K^2_S$ as complete intersections in projective
bundles over elliptic curves.

Using Atiyah's classification of vector bundles over elliptic curves, and the
representation theory of the Heisenberg groups, it was possible to
obtain the following result

\

{\bf  Theorem (\cite{cc1}, \cite{cc2})} {\it The Albanese fibre of a 
minimal surface
with $K^2_S = 3, p_g= q=1$  has genus  $g=2$ or $g=3$.

Surfaces with $g=3$, whose canonical model is a divisor $\Sigma$
in the third symmetric product $B^{(3)}$, yield a unirational 
connected component of
the moduli space of dimension $5$.

A surface for which $g=2 $ is a double cover of the symmetric product
$B^{(2)}$  branched in a divisor $\Delta$ belonging to a fixed
algebraic system of 6-sections, and with two 4-uple singular points 
on the same fibre.
These surfaces exist. }

\

\cite{cc1} asked whether the second family is irreducible (and
conjectured somewhat overhastily that the answer should be positive).
A main application of our structure theorem for genus $g=2$ is the complete
classification of the moduli space of the above surfaces. This
classification, while in particular  answering  the above
question in the negative, shows however the existence of an irreducible
{\bf main stream family } of  surfaces with $K^2_S
= 3, p_g= q=1, g=2.$

\

  {\bf Theorem \ref{Kquadro3gen2} } {\it The moduli space $\M$ of minimal
surfaces of general type with $K^2_S = 3, p_g= q=1$  has exactly
four connected components, all (irreducible) unirational of dimension
$5$.}

\begin{oss}
\cite{isogenous} shows, as a consequence of Seiberg Witten theory, 
that the genus $g$ of the Albanese fibre is a differentiable invariant of the 
underlying complex surface.

Thus restricting ourselves  to the case $g=2$, we may ask about the 
existence of surfaces with $p_g= q=1, K^2_S > 3$ and $g=2.$ In this
case, by a result of Xiao  (\cite{xiao})
one has  the inequality $ K^2_S \leq 6$, and \cite{bidouble}
proved the existence of surfaces with  $p_g= q=1,  K^2_S = 4,5$ and $g=2.$
\end{oss}

We wish to propose  the following questions as  further applications
of our methods.

{\bf Problem I:} Do there exist surfaces with $p_g= q=1, K^2_S = 6$
and $g=2$?

{\bf Problem II:} How many connected and irreducible components do the
moduli spaces corresponding to surfaces with $p_g=q=1$, $K^2_S=4,5, 6$
and $g=2$ have?

We  now  introduce our structure theorems for
fibrations of genus 2 and 3. We hope to treat the case $g=4$
in a sequel to the present article (motivation here
came from the classification problem of numerical Godeaux surfaces,
cf. \cite{reidtor} and \cite{cp}).

Since the canonical map of a curve of genus 2 is  a double cover of $\PP^1$
branched in 6 points, one of the first geometric approaches, developped
by Horikawa, was to study  genus $2$ fibrations $f : S \ra B$
through a finite double cover $ \phi :Y \ra  \PP $ of a $\PP^1$ bundle
$\PP$ over $B$, where $Y$ is birational to $S$.

There are two drawbacks here: first, $Y$ and $\PP$ are not uniquely determined,
second,  the branch curve $\Delta$ has several {\bf nonsimple} 
singularities, which
are usually several pairs of 4-tuple points or of (3,3) points occurring on
the same fibre.

Our approach uses the geometry of the bicanonical map of a 1-connected divisor
of genus 2, which is a morphism generically of degree 2 onto a plane conic $Q$
which may be reducible or nonreduced.

In this way we obtain a unique birational model $X$ of $S$,
admitting a finite double cover $\psi : X \ra \CCC$,
where

(i)  $\CCC$ is   a {\bf conic bundle}  over $B$

(ii) the branch curve $\Delta_{\A}$ has only {\bf simple}
singularities

(iii) $X$ is the {\bf relative canonical model} of $f$, and is 
obtained contracting
the (-2)-curves contained in the fibres to singularities which are then
Rational Double Points.

In order to better explain the meaning of (iii) above, and since
our approach ultimately provides
purely algebraic structure theorems, we need to recall the algebraic 
definition of
  the
{\bf relative canonical algebra}  $\R(f)$ of the fibration (its local structure
  was investigated for $ g \leq 3$ by Mendes Lopes in \cite{mml} and its
importance for
general $g$ was stressed in
\cite{reid}).

  To state our two main results we need some terminology:
we consider the canonical ring of any fibre $F_t$ of $ f : S \ra B$
$$  \R(F_t) :=  \oplus_{n=0}^{\infty} V_n(t) := \oplus_{n=0}^{\infty} H^0 (F_t,
\hol_{F_t} ( n K_{F_t}) ),$$ 
where $K$ denotes the canonical divisor. These
vector spaces fit together, yielding vector bundles $V_n$   and the relative
canonical algebra $ \R (f)$ on $B$, defined as follows:
$$ \R (f) : = \oplus_{n=0}^{\infty} V_n : = \oplus_{n=0}^{\infty}
f_*(\hol_S ( n ( K_S - f^* (K_B)))).$$

Write $\sigma_n : Sym^n (V_1) \ra V_n $ for the multiplication
map.   The sheaf $ \T_2 : = Coker (\sigma_2)$
introduced in \cite{cc1} also plays an important role, as does the
  class $\xi$ of the extension
$$ 0 \rightarrow Sym^2(V_1) \rightarrow V_2 \rightarrow \T_2
\rightarrow 0;$$
we show here that $ \T_2$ is isomorphic to the structure
sheaf $\hol_{\tau}$ of an effective divisor $\tau$ on the curve $B$.

The geometry behind the above exact sequence is that $\sigma_2$
determines a (rational) relative quadratic Veronese map $\PP(V_1) 
\dashrightarrow \PP(V_2)$ whose image is the conic bundle $\CCC$
mentioned above. We denote 
by $\A$ the relative anticanonical algebra of the conic bundle $\CCC$,
and we show that $\R$ is a locally free $\A$-module, $ \R \cong \A 
\oplus (\A[-3] \otimes (\det V_1 \otimes \hol_B(\tau))$. 

The final datum, denoted by $w$, contains
  the most geometric meaning: when $g=2$ it determines
the branch divisor $\Delta_{\A}$ (intersecting the  fibres
of the conic bundle $\CCC$ in $6$ points), when $ g=3$ it
determines a divisor $\Sigma$ in the $\PP^2$ bundle $\PP (V_1)$
intersecting a general fibre $\cong \PP^2$ in the quartic curve
$\Sigma_t$ which is the canonical image of the nonhyperelliptic
fibre $F_t$.

Algebraically, when $g=2$, $w \in  \PP ( {\rm Hom}((\det V_1 \otimes 
\hol_B(\tau))^2,
\A_6) )$ gives the multiplication map
$ \big(\A [-3] \otimes
(\det V_1 \otimes 
\hol_B(\tau))\big)^{\oplus 2} \ra \A$, whence it determines the  algebra 
structure of
$\R$.

Denote by $b$ the genus of $B$, and define the condition of admissibility
of a 5-tuple $(B,V_1,\tau,\xi,w)$ as above as the open condition
which guarantees that the singularities of the  relative canonical model $X : =
{\mathbf{Proj}} (\R)$  are Rational Double Points.
We can now give  the exact statements of our two structure
theorems.

\

{\bf Theorem \ref{gen2}}
{\it Let $f$ be a relatively minimal genus $2$ fibration. Then its
{\bf associated 5-tuple} $(B,V_1,\tau,\xi,w)$ is admissible.

Viceversa, every admissible genus two $5-$tuple determines a
sheaf of algebras $ \R \cong \A \oplus ( \A [-3] \otimes \det(V_1) 
\otimes \hol_B
(\tau))$ over
$B$ whose relative projective spectrum $X$ is the relative canonical model
of a relatively minimal genus $2$ fibration $f:S \rightarrow B$
having the above as associated 5-tuple. Moreover, the
surface $S$ has
the following invariants:
$$\begin{array}{lll} \chi(\hol_S)&=&\deg(V_1)+(b-1)\\ K_S^2&=&2\deg V_1
+ \deg \tau + 8(b-1).
\end{array}$$

We obtain thus a bijection between (isomorphism classes of)
  relatively minimal genus $2$
fibrations and  (isomorphism classes of) associated
$5-$tuples are isomorphic, which is functorial in the sense that to a flat family of
fibrations corresponds a flat family of 5-tuples.}

\

{\bf Theorem \ref{strucgen3} }
{\it
  Let $f$ be a relatively minimal genus $3$ nonhyperelliptic
fibration such that every fibre is $2-$connected.
Then its
{\bf associated 5-tuple} $(B,V_1,\tau,\xi,w)$ is admissible.

Viceversa, every admissible genus three 5-tuple
($B$, $V_1$, $\tau$, $\xi$, $w$)
is the associated 5-tuple of a unique  genus
$3$ nonhyperelliptic fibration
$f:S \rightarrow B$ with the property that every fibre is $2$-connected
and with invariants $\chi (\hol_S) =\deg V_1 +
2(b-1)$, $K^2_S=3\deg V_1 + \deg \tau +16 (b-1)$.
As in the case of genus 2, the bijection thus obtained is  functorial.
}

\

In the somewhat shorter last section we give an  application  of
the structure theorem for
$g=3$ in a case where the base curve $B$ is $\PP^1$. This case
is easier than the case where the genus $b$ is higher since
by Grothendieck's theorem every vector bundle on $\PP^1$ is a direct sum
of line bundles. Our theorem, establishing a new result, namely
the existence
of nonhyperelliptic genus $3$ fibrations for some subvarieties of the
moduli spaces of surfaces with
$ p_g =3, q=0$ and $ K^2_S = 2,3,4,5$, should be viewed as a guidebook
to the use of our structure theorems for the case where
$ B = \PP^1$.

As already mentioned, we plan in a sequel to this paper to describe the case of
hyperelliptic fibrations of genus $g=3$ and the case of nonhyperelliptic
fibrations of genus
$g=4$,  giving applications to other  questions of
surface theory.

We also hope that  our present results may be found useful in
  developing the arithmetic theory of curves of genus $g=2$,
resp. of genus $g=3$.

\section{Generalities on fibrations of surfaces to curves}

Throughout this paper $f: S \rightarrow B$
will be a relatively minimal fibration
of a projective complex
surface $S$ onto a smooth projective complex curve $B$ of genus $b$.

This means that  $f$ has connected fibres $F$, and that no fibre contains an
exceptional curve of the first kind.

We denote by $g$ the  arithmetic genus of $F$, and
we shall assume $g \geq 2$.

For every $p \in B$, we shall denote by $F_p$ the fibre of $f$ over
$p$, i.e., the divisor $ f^*(p)$. The canonical divisor $K_{S}$
restricts on each fibre to the dualizing sheaf $\omega_{F_p}$; recall
that, $f$ being a fibration, $\forall p \in B$, $h^0(\hol_{F_p})=1$
(see  \cite{bpv}, lemma III 11.1) and $h^0(\omega_{F_p})=$ genus
$(F_p)=:g$.

Assume that $S'$ is a smooth projective surface  and that $f':S'
\dashrightarrow B$ is a rational map onto a smooth  curve $B$.

$f'$ is necessarily a morphism if $ b$:= genus $(B)$ is strictly
positive: when $b$ is equal to zero, we let $ \beta: S \rightarrow S'$
be a minimal sequence of blow ups which yields a morphism $f: S
\rightarrow B$ and the following commutative diagram
\[
      \begin{array}{rcl}
       \label{diag}
     S &  & \\ \beta \downarrow & \searrow \! \! \!
\stackrel{\stackrel{\scriptstyle f}{}}{} & \\ S' &
\stackrel{f'}{\dashrightarrow} & B.
      \end{array}
\]

If $S'$ is a minimal surface of general type, $f'$ is "almost always" 
a morphism,
we have in fact (cf. Kodaira's lemma in \cite{horp} and \cite{xiao}, 
prop. 4.1))

\begin{lem}\label{minimalgen2}
A relatively minimal genus $2$ fibration on a non minimal surface $S$
of general type occurs only for the canonical pencil of a minimal
surface $S'$ with $K_{S'}^2=1, p_g (S') =2.$
\end{lem}

\

Although less explicit, the next result also gives rise to a finite
number of families.

\begin{lem}\label{minimalgeng}
A relatively minimal genus $g$ fibration on a non minimal surface $S$
of general type can only occur if its minimal model
  $S'$ has $$K_{S'}^2 \leq ( 2g-3)^2.$$
\end{lem}
\Proof
We have $K_S = \beta^* K_{S'} + \sum_{i=1,..r} E_i$, where each $E_i$
is an exceptional divisor of the first kind (i.e., the total transform
of the maximal ideal of a smooth point).  We have $ 2 g - 2=
K_S F = F
\beta^* K_{S'} + \sum_{i=1,..r} F E_i$.  Let $M$ be the corresponding
pencil on $S'$: thus $F \beta^* K_{S'} = K_{S'} M \geq 1$. On the other hand,
$F E_i \geq 1$, else $E_i$ is vertical and the fibration is not
relatively minimal.  Set $\beta^* (M) = F -  \sum_{i=1,..r} m_i E_i$,
so that $ m : = \sum_{i=1,..r} F E_i = \sum_{i=1,..r} m_i \leq 2g -
3$.

We obtain $ M^2 \leq m^2$, and by the index theorem follows then
$ K_{S'}^2 \cdot M^2 \leq (K_{S'} M)^2 $, whence
$ K_{S'}^2  \leq (2g -2 -m)^2  \leq (2g -3)^2$.
\qed

The canonical maps of the fibres can be combined into the
{\em relative canonical map}(see e.g. \cite{horp}), given concretely 
as follows: let $L$ be a divisor on $B$ such that $H^0(S,\hol_S(K_S +
f^{*}L)) \rightarrow H^0(F_p,\omega_{F_p})$ is surjective $\forall p
\in B$, and consider the rational map
$$h:S \rightarrow \PP({H^0(S,\hol_S(K_S + f^{*}L)}) \times B,$$
induced by the linear system $|\hol_S(K_S + f^{*}L)|$ and
the projection $f$.

$h$ is a birational map unless every fibre is hyperelliptic: in this
case $h$ is a double cover of a surface $Y$ ruled over $B$.
\begin{ex}
\label{easy}
If $B=\PP^1$ and $p_g(S)=q(S)=0$, then,  setting  $L : =F$, we
obtain $h : S \rightarrow \PP^{g-1} \times \PP^1$, and every fibre
of $f$ is mapped to $\PP^{g-1}$ via its canonical map.
\end{ex}
\Proof For every $p \in \PP^1$, the
exact sequence
\[
0 \rightarrow \omega_S \rightarrow \hol_S(K_S +F)
\rightarrow \omega_{F_p} \rightarrow 0
\]
yields  an isomorphism $H^0(\hol_S(K_S +F)) \rightarrow
H^0(\omega_{F_p})$.
\qed

The previous example is particularly relevant in the case of numerical
Godeaux surfaces, to which we are interested (cf. Theorem
\ref{godeaux}).

Let us first collect some known results on canonical maps of
Gorenstein curves ( cf. \cite{cf}, \cite{4names}, \cite{cp},
\cite{mml}).

These are based on Franchetta's definition (cf. \cite{fra1},  \cite{fra2})
\begin{df}
An effective divisor $D$ on a smooth algebraic surface $S$ is said to
be $k$-connected if, whenever we write $D= A+B$ as a sum of effective
    divisors $A,B > 0$, then $ A \cdot B \geq k$.
\end{df}
\begin{lem}
\label{2-conn}
Let $F$ be a $2$-connected curve of genus $ g \geq 2$:
then the canonical map of $F$ is
a morphism.
\end{lem}
\begin{lem}
\label{3-conngen3}
Let $C$ be a curve of genus $ g \geq 3$ , let $\omega$ be the
dualizing sheaf of
$C$.  Then $\varphi_\omega$ embeds $C$
$\Leftrightarrow$
$C$ is $3$-connected and $C$ is not honestly hyperelliptic (i.e., a finite
double cover of $\PP^1$ induced by the canonical morphism).
\end{lem}
\begin{lem}\label{2-conngen3}
Let $C$ be  $2$-connected of genus $ g = 3$ with $\omega$ ample
(e.g., if $C$ is a fibre of
the relative canonical model of a fibration
of a surface of general type ):
then either  $\varphi_\omega$ embeds $C$ as a plane
quartic, or $\varphi_\omega$ is a finite double cover of a plane conic,
and more precisely $C$ is a complete intersection of type $(2,4)$
in the weighted projective space $\PP (1,1,1,2)$.
\end{lem}
\begin{lem}
\label{3-conn}
Let $C$ be a $3$-connected genus $4$ Gorenstein curve, which is not honestly
hyperelliptic: then $\varphi_\omega$ embeds
$C$ as a complete intersection of type $(2,3)$.
\end{lem}
Let us now recall an important definition, of the relative canonical
algebra $\R(f)$ of the fibration $f$ (cf. \cite{cc1}, \cite{reid})
\begin{df}
Consider the relative dualizing sheaf
$$\omega_{S|B}:=\hol_S( K_S-f^*K_B).$$
Its self-intersection defines the number
  $$ K^2_{S|B}: = (\omega_{S|B})^2=K_S^2-8(b-1)(g-1)$$ and we further define
$\chi(S|B) : =\chi(\hol_S)-(b-1)(g-1)$.
The relative canonical algebra $\R(f)$ is
the graded algebra $$\R(f) := \oplus_0^{\infty} V_n,$$ where $V_n$ is 
the vector
bundle on
$B$ given by $ V_n : = f_*(\omega_{S|B}^{\otimes n})$.
\end{df}
\begin{oss}\label{fujita}
Since we have a fibration, $V_0 = \hol_B$ , $V_1$ has rank $g$,
$V_n$ has rank $(2n-1)(g-1)$.

By a theorem of Fujita (\cite{fu1},
\cite{fu2}) $V_n$ is semipositive,  meaning that every
rank $1$ locally free quotient of it has non negative degree,
and $V_1$ is a direct sum $\hol_B^{q(S)-b} \oplus V_1^a \oplus
V_1^0$ where $V_1^a$ is ample and $V_1^0$ is a direct sum of stable
degree 0 bundles $\E_i$ with $H^0(\E_i)=0$;
if rank $ \E_i=1$, then $\E_i$ is a line
bundle associated to a torsion divisor in $Pic^0(B)$ (
cf. \cite{zuc}).

Notice that, by relative duality, $R^1f_* \omega_{S|B}=\hol_B$, while
$R^1f_* \omega_{S|B}^{\otimes n}=0$ for $n \geq 2$ because the fibration is
assumed to be relatively minimal.

It then follows  by Riemann-Roch that for $n \geq 1$
$$
\chi(V_n)
=\chi(\omega^{\otimes n}_{S|B})=
\frac{1}{2}n(n-1)K^2_{S|B}+ rank(V_n) \cdot \chi(\hol_B)+
\chi(S |B).$$
\end{oss}
\begin{oss}\label{fujitadeg}
Fujita's theorem shows that
$$
deg(V_n)= \frac{1}{2}n(n-1)K^2_{S|B}+
    \chi(S |B) \geq 0.$$
The Arakelov inequality
\[ K^2_{S|B} = [K_S^2-8(g-1)(b-1)] \geq 0\]
   follows  as a
corollary, together with the inequality
\[ \chi(S |B): = \{ \chi(\hol_S)-(g-1)(b-1)\} \geq 0 .\]
Moreover (cf. \cite{e-v}) $V_n$ is  ample for $ n \geq 2$ if
the fibration is
  a nonconstant moduli fibration (i.e.,  the smooth fibres
are not all isomorphic), and in this case
\[ K^2_{S|B} = [K_S^2-8(g-1)(b-1)]> 0.\]
Moreover, in the inequality
\[\chi(\hol_S)-(g-1)(b-1) \geq 0,\] equality holds if and only
if $f$ is a holomorphic bundle (by Noether's formula,
cf. \cite{beau}).

Finally, we have the inequality $q(S) \leq b + g$, equality holding
if and only if $f$ is a product fibration $ F \times B \rightarrow B$
(see \cite{beau}).
\end{oss}

\section{Invariants of the relative canonical algebra }
Consider now the relative canonical algebra $$ {\R} (f)= \oplus_0^{\infty}
V_n= \oplus_0^{\infty} f_*(\omega^{\otimes n}_{S|B}).$$
\begin{df}
Denote by $\mu_{n,m}: V_n \otimes V_m \rightarrow
V_{n+m}$, respectively by
\[ \sigma_n:
S^n(V_1) : = Sym^n(V_1) =S^n(f_*(\omega_{S|B}))
\stackrel{\sigma_n}{\longrightarrow}V_n= f_*(\omega^{\otimes n}_{S|B}),
\]
the homomorphisms induced by multiplication.

We define  $\L_n : =\ker{\sigma_n}$ and $\T_n := \hbox{ \rm
coker}\ \sigma_n$.
\end{df}
\begin{oss}\label{torsione}
By Noether's theorem on canonical curves, $\T_n$ is a torsion sheaf if
the general fibre of $f$ is nonhyperelliptic, whereas more generally,
if the general fibre is nonhyperelliptic, coker $\mu_{n,m}$ is a
torsion sheaf as long as $n,m \geq 1$.

A more precise result was recently proved by K. Konno and
M. Franciosi (\cite{kon}, \cite{fran}): if every fibre is
$1$-connected (equivalently,
there is no multiple fibre), coker $\mu_{n,m}$ = 0 if $ g\geq 2, n,m \geq 2,
max\{n,m\} \geq 3$.

When there are no multiple fibres, the relative canonical
algebra is generated by elements of degree $ \leq 3$, and in
general it is generated by elements of degree $\leq 4$.
\end{oss}
The previous remark shows that the two cases
\begin{itemize}
\item
I) The general fibre is nonhyperelliptic
\item
II) All the fibres are hyperelliptic.
\end{itemize}
should be treated separately.

In the hyperelliptic case, one has the following useful
method, of splitting the relative canonical algebra into the
invariant, resp. anti-invariant part.
Assume in fact  that a general fibre is
hyperelliptic: then there is a birational involution $\sigma$ on $S$,
whence also on $S'$. Assuming that $S$ is not birationally ruled,
$\sigma$ acts biregularly on $S'$. Moreover, since $\sigma$ preserves
the rational map $f'$, it preserves its indeterminacy locus, therefore
the minimal sequence of blow-ups turning $f'$ into a morphism is
$\sigma$-equivariant, hence one concludes that $\sigma$ acts biregularly
on the fibration $ f: S \rightarrow B$ (and trivially on the base
$B$).

Therefore, each open set $U = f^{-1} U'$ is $\sigma$-invariant and
$\sigma$ acts linearly on the space of sections
$\hol_{S}(U,\omega_{S/B}^n)$, which splits as the direct sum of the
$(+1)$-eigenspace and the $(-1)$-eigenspace.

Accordingly, we get direct sums:
$$V_n= V_n^+ \oplus V_n^- = f_*(\omega^{\otimes n}_{S|B})^+ \oplus
f_*(\omega^{\otimes n}_{S|B})^-$$ and we can split the relative canonical
algebra as:
$$ \R (f)= \R (f)^+ \oplus \R (f)^-,$$ where we now observe that $ \R
(f)^+$ is a subalgebra and $\R (f)^-$ is an $ \R (f)^+$-module.
\begin{oss}\label{ranks}
Let the general fibre of the fibration $f$ be hyperelliptic:
then $V_1= V_1^- $, and the cokernels $\T_n $ will be bigger than in the non
hyperelliptic case.

We have in fact the following table for the ranks of $V_n^+$
and $V_n^-$ (for $n \geq 2$):
\begin{center}
\begin{tabular}{|l|l|l|}
\hline n & rank $V_n^+$ & rank $V_n^-$ \\ \hline even &$n(g-1)+1$ &
$(n-1)(g-1)-1$\\ odd &$(n-1)(g-1)-1$ & $n(g-1)+1$\\ \hline
\end{tabular}
\end{center}
Then the sheaf maps
\[
S^n(V_1) =S^n(f_*(\omega_{S|B}))
\stackrel{\sigma_n}{\longrightarrow}V_n= f_*(\omega^{\otimes n}_{S|B}),
\]
are injective iff $ g = 2$ and their image lies in $V_n^+$ for $n$ even and in
$V_n^-$ for $n$ odd.

So, if we define $ \T_n^+=$ coker $S^n(f_*(\omega_{S|B})) \rightarrow
V_n^+$ for $n$ even, and $\T_n^-=$coker $S^n(f_*(\omega_{S|B}))
\rightarrow V_n^-$ for $n$ odd, the decomposition in invariant and
anti-invariant part of the cokernels $\T_n$ is $\T_n=\T_n^+ \oplus
V_n^-$ for $n$ even, $\T_n=\T_n^- \oplus V_n^+$ for $n$ odd,
and the sheaves $\T_n^{\pm}$  are torsion sheaves.
\end{oss}
We end this section by observing that
$$
\deg (Sym^n(V_1)) = {n+g-1 \choose g} \deg (V_1)  ={n+g-1 \choose g}  \chi
(S | B).$$

\section{Genus 2 fibrations: the structure theorem}
Before describing  the building data of the
relative canonical algebra of a genus $2$ fibration, it is
convenient to explain the underlying geometry (cf. also
\cite{cc1},\cite{cc2}).

Let $f:S \rightarrow B$ be a genus $2$ fibration.  The rank $2$ vector
bundle $V_1:=f_* \omega_{S|B}$ induces a natural factorization of $f$
as $\pi \circ \varphi$, where $\varphi: S \dashrightarrow
\PP : = \PP(V_1) =\Proj(Sym(V_1))$ is a rational map of degree $2$ , and $\pi:
\PP(V_1) \rightarrow B$ is the natural projection.

The indeterminacy locus of $\varphi$ is contained in the fibres of $f$
that are not 2-connected, i.e.,  split as $F_p = \E_1 + \E_2$ with
$\E_1  \E_2 = 1$. Then $\E_i ^2 = -1$, $\E_i$ has arithmetic genus
$1$ and is called an {\em elliptic cycle}. We will see that these
fibres are exactly the inverse images of the points in $Supp (\T_2)$.

The typical example is given by a fibre consisting of two smooth
elliptic curves $\E_1, \E_2$ meeting transversally
in a point $P'$. The blowup of  $P'$ maps isomorphically to
the fibre $F''$ of $\PP $ over the point $ P \in B$, while the
elliptic curves $\E_1, \E_2$ are contracted to two
distinct points of the fibre $F''$.
The resolution $\tilde{\varphi}$ of $\varphi$ is the composition of
the contraction of $\E_1, \E_2$ to two simple
$-2$-elliptic singularities, with a finite double cover where the
branch curve $\Delta$ in $\PP$ contains the fibre and has two
distinct $4$-tuple
points on it. More complicated fibres containing elliptic tails can
produce different
configurations of singularities of the branching divisor of
$\varphi$:  complete
lists are given by  Ogg (  \cite{ogg}) and by Horikawa
  (\cite{horp}).
\begin{lem}\label{torsgenus2}
Let $f:S \rightarrow B$ be a genus $2$ fibration.  Then
\begin{enumerate}
\item $\T_2$ is isomorphic to the structure sheaf of an effective 
divisor $\tau \in
    Div_{\geq 0}(B)$,  supported on the points of $B$
    corresponding to the fibres of $f$ which are not 2-connected;
\item $\tau$ determines all the (torsion) sheaves $\T_n$ as follows:
$$\T_{2n}^+ \cong \hol_{n\tau} \oplus (\oplus_{i=1}^{n-1} \hol_{i
     \tau})^{\oplus 2} \ \ \ \ \T_{2n+1}^- \cong (\oplus_{i=1}^n \hol_{i
     \tau})^{\oplus 2};
$$ In particular $\deg \T_{2n}^+ = n^2 \deg \tau$ and $\deg
\T_{2n+1}^- = n(n+1) \deg \tau$.
\end{enumerate}
\end{lem}
\Proof As shown in the thesis of M.Mendes Lopes (thm. 3.7 of \cite{mml}, page
53), there are two possibilities for the canonical ring of a genus $2$
fibre:
\begin{itemize}
\item either the curve is honestly hyperelliptic, i.e., the graded
ring is isomorphic
to
$$ \C [x_0,x_1,z]/ (z^2 -g_6(x_0,x_1)),$$ where $ \deg x_0= \deg x_1=
1$, $\deg z=3$, $\deg g_6 = 6 $,
\item or the fibre is not $2$-connected and the ring is isomorphic to
$$ \C [x_0,x_1,y,z]/(Q_2, Q_6),$$ where $\deg x_0= \deg x_1= 1$, $\deg
y = 2$, $\deg z=3$ and

$Q_2: = x_0^2 - \lambda x_0 x_1$

$Q_6: = z^2 - y^3 - x_1^2 (\alpha_0 y^2 + \alpha_1 x_1^4) $.
\end{itemize}
The first case is the one where the fibre is $2$-connected.

By remark \ref{ranks} $ x_0, x_1$ are anti-invariant sections,  $y,z$ are
invariant and the sheaves
$\T_{2n}^+$, $\T_{2n+1}^-$ are zero away from the points $P$ whose fibres
are not $2$-connected.

$V_2^+$is locally generated by $ x_0^2,x_0 x_1,x_1^2,y $.

By flatness, if $t$ is a uniformizing parameter for $\hol_{B,P}$,
we can lift the relation
$Q_2: = x_0^2 - \lambda x_0 x_1$ to
$$(x_0^2 - \lambda x_0 x_1) + t \mu(t) y + t [x_0^2 \psi_0 (t) + x_1^2 \psi_1
(t) + x_0 x_1 \phi (t)] : = Q_2 + t \mu(t) y + t R(x,t) .$$

$\mu (t)$ is not identically $0$ since  $x_0$ and $ x_1$ are algebraically
independent on a general fibre.
Therefore, by a holomorphic change of coordinates in $B$, we may
assume $\mu(t) = t^{s-1}$
for a suitable positive integer $s$: we will call this integer the
''multiplicity'' of our special point.  The above relation
  shows that the stalk of $\T_2 $ at a special point $P$ is the
principal module $\hol_{B,P}/ (t^s)$ generated by the class of $y$.

We can also choose a lifting $Q_6(t)$ of $Q_6$ of the form
$$Q_6(t): = z^2 -Q'_6(x_0, x_1,y, t):$$
since $Q_6$ is invariant, the lifting $Q_6(t)$ must be invariant
too, otherwise  a nontrivial antiinvariant relation would imply
that $z$ is identically zero, absurd.

By flatness these are all the relations of the stalk of $\R$ at $P$; we
leave to the reader the straightforward computation showing that
$(\T_{2n+1}^-)_P$ equals $\oplus_{i=1}^n \left( \hol_{B,P}/(t^{is})
\oplus \hol_{B,p}/(t^{is}) \right)$ with minimal ordered system of
generators
$$\{x_0x_1^{2n-2}y,x_1^{2n-1}y,\ldots,x_0y^n,x_1y^n\},$$ and
$(\T_{2n}^+)_P$ equals $ \left( \oplus_{i=1}^{n-1} \left(
\hol_{B,p}/(t^{is}) \oplus \hol_{B,p}/(t^{is}) \right)\right) \oplus
\hol_{B,p}/(t^{ns})$ with minimal ordered system of generators
$$\{x_0x_1^{2n-3}y,x_1^{2n-2}y,\ldots,x_0x_1y^{n-1},
x_1^2y^{n-1},y^n\}.$$
\qed

We now introduce  a sheaf of graded algebras, whose $\Proj$ yields
a conic bundle $\CCC$
over $B$ admitting the relative canonical model $X$ as a finite double cover.
\begin{df}\label{A}
Let  $\A$ be the graded subalgebra of $\R$ generated by $V_1$ and
$V_2$.

$\A_n$  denotes its graded part of degree $n$, and
write accordingly
$$\A= \A_{even}\oplus \A_{odd} := (\oplus_{n=0}^{\infty} \A_{2n})
\oplus (\oplus_{n=0}^{\infty} \A_{2n+1}).$$ 
We decompose similarly  $\R=\R_{even} \oplus \R_{odd}$.
\end{df}
\begin{lem}\label{splittingR}
$\R$ is isomorphic to $\A \oplus (\A[-3] \otimes V_3^+)$ as a graded
$\A$-module; moreover $\A_{even}$ is the invariant part of $\R_{even}$
and $\A_{odd}$ is the antiinvariant part of $\R_{odd}$.
\end{lem}
\Proof In the proof of lemma \ref{torsgenus2} we wrote the stalk of
$\R$ at a special point $P$ as an $\hol_{B,P}$-algebra generated by
$x_0, x_1,y$ and $z$, where the $x_i$'s are antiinvariant of degree
$1$, $y$ and $z$ are invariant with respective degrees $2$ and $3$.
We achieve a unified treatment of both cases writing the
canonical ring of
a honestly hyperelliptic fibre as
$$ \C[x_0,x_1,y,z]/(y,z^2-g_6) : = \C[x_0,x_1,y,z]/ (Q_2, Q_6).
$$

In both case the stalk of $\A$ is the subalgebra generated by
$x_0,x_1$ and $y$: therefore $\A_{even}$ is invariant and $\A_{odd}$
is antiinvariant.

Locally on $B$ we have $\R \cong  \hol_B [x_0,x_1,y,z]/ (Q_2 (t),
Q_6 (t)) $, and  since $ Q_6(t) = z^2-Q_6'(x_0,x_1,y,t)$,
it follows easily that locally

(1) $\A \cong  \hol_B [x_0,x_1,y]/ (Q_2 (t))$

(2)  $\R \cong \A \oplus z \A$.

Since $z$ is a local generator of $V_3^+$, both assertions follow.
\qed

In the next lemma we give a more explicit description of $\A$,
showing how we can construct $\A$ starting from $\sigma_2$.
\begin{lem}\label{structA}
There are exact sequences
\begin{equation}\label{exseqAeven}
0 \rightarrow (\det V_1)^2 \otimes S^{n-2}(V_2)
\stackrel{i_n}{\rightarrow} S^n(V_2) \rightarrow \A_{2n} \rightarrow 0
\ \ \ \forall n \geq 2;
\end{equation}
\begin{equation}\label{exseqAodd}
    V_1 \otimes (\det V_1) \otimes \A_{2n-2} \stackrel{j_n}{\rightarrow}
V_1 \otimes \A_{2n} \rightarrow \A_{2n+1} \rightarrow 0 \ \ \ \forall
n \geq 1;
\end{equation}
where
$$i_n((x_0 \wedge x_1)^{\otimes 2} \otimes
q):=(\sigma_2(x_0^2)\sigma_2(x_1^2) -\sigma_2(x_0x_1)^2)q,$$
$$j_n(l \otimes (x_0 \wedge x_1) \otimes q):=x_0 \otimes
(\sigma_2(x_1l)q)-x_1 \otimes (\sigma_2(x_0l)q).$$
\end{lem}
\Proof The above maps $S^n(V_2) \rightarrow \A_{2n}$ and $V_1 \otimes
\A_{2n} \rightarrow \A_{2n+1}$, induced by the
ring structure of $\A$,  are surjective because  $\A$
is generated in degree $\leq 2$ by definition. Since
$\R_n$ and $\A_n$  are locally free, the respective kernels  are
locally free sheaves on $B$.

Both sequences are complexes by virtue of associativity and
commutativity of multiplication in $\R$.

To verify exactness, we first consider the
stalk at a special point
$P \in Supp (\T)$. The sheaf homomorphisms $i_n$, $j_n$ and $\sigma_2$
induce linear maps $i_{n,P}$, $j_{n,P}$, $\sigma_{2,P}$ on the fibres
over $P$ of the corresponding vector bundles, and ring homomorphisms
$i_{n,(P)}$, $j_{n,(P)}$ , $\sigma_{2,(P)}$ on the stalks of the
corresponding sheaves.

In the proof of lemma \ref{torsgenus2} we wrote the relation in degree
$2$ of ${\A}$ as $t^sy-Q'(x_0,x_1,t)$ with
$Q'(x_0,x_1,0)=-x_0^2+\lambda x_0 x_1$.

Therefore $u_0:=y,u_1:= \sigma_{2,(P)}(x_0x_1)$ and
$u_2:=\sigma_{2,(P)}(x_1^2)$ are a basis of the fibre of the stalk
(and by restriction of the fibre) of $V_2$ at $P$, and
$\sigma_{2,(P)}(x_0^2)=-t^su_0+\lambda u_1 +t\cdots$.

In this basis $i_{n,P}((x_0 \wedge x_1)^2 \otimes q)= (\lambda u_1u_2
-u_1^2)q$, thus $i_{n,P}$  is clearly injective.
At a general point we choose
  $$u_0:=\sigma_{2,(P)}(x_0^2), \
u_1:=\sigma_{2,(P)}(x_0x_1)\ {\rm and} \ \  u_2:=\sigma_{2,(P)}(x_1^2)$$
as basis and since
$i_{n,P}((x_0 \wedge x_1)^2 \otimes q)= ( u_0u_2 -u_1^2)q$ we derive the same
conclusion.

Since $i_n$ injects $(\det V_1)^2 \otimes S^{n-2}(V_2)$
in $S^n(V_2)$ as a saturated subbundle, the exactness of the sequence
(\ref{exseqAeven}) follows then from the equality
\begin{multline*}
{\rm rank} ((\det V_1)^2
\otimes S^{n-2}(V_2)) + {\rm rank} (\A_{2n})={n \choose 2}+2n+1= {\rm rank}
(S^n(V_2)).
\end{multline*}

To show that the image of $j_n$ contains the kernel of the
projection on $\A_{2n+1}$, which is locally free by our previous remark,
   it is enough to work on the fibres of the associated vector bundles.
We identify the fibre of $\A_{k}$ with a subspace of the canonical ring
of the fibre curve.
We must then  show that, given $p_0,p_1$ of
even degree with $x_0 p_1=x_1 p_0$, there exists a $q$ with $
p_i=x_iq$ ($i=0,1$).

This is straightforward for the fibre of a general point, because this
subring is the ring $\C[x_0,x_1]$.

On a special fibre our subring is the ring $ \C[x_0,x_1, y]/(x_0^2-\lambda
x_0x_1)  = (\C[x_0,x_1]/(x_0^2-\lambda
x_0x_1))[y] $. We can thus replace $p_0,p_1$ by their coefficients in
   $\C[x_0,x_1]/(x_0^2-\lambda x_0x_1)$.

In this ring, if $k \geq 2$ is the degree of the polynomials $p_i$, we can
  find uniquely determined constants $a_i,b_i$ with
$p_0=a_0x_1^{k-1}x_0+b_0x_1^{k}$, $p_1=a_1x_1^{k-1}(x_0-\lambda
x_1)+b_1x_1^{k}$.
With this expression of the polynomials the condition $x_0p_1-x_1p_0$
becomes $b_0=a_0-b_1=0$ and it suffices to choose
$q=a_0x_1^{k-1} +a_1x_1^{k-2}(x_0-\lambda x_1)$.
\qed

\begin{oss}\label{sigma2determinesA}
The above exact sequences describe the subalgebra $\A_{even}$
as a quotient algebra of $Sym(V_2)$ and $\A_{odd}$
as an $\A_{even}$ module. The multiplication map
$\A_{odd} \times \A_{odd} \ra \A_{even}$
is induced by $$ \mu_{1,1} : V_1 \otimes V_1 \ra S^2 (V_1) 
\stackrel{\sigma_2}{\ra}
V_2.$$

Thus ${\A}$ is completely determined as an $\hol_B$-algebra by
the map $\sigma_2:S^2(V_1) \rightarrow V_2$.
\end{oss}

The structure of quotient algebra of $Sym(V_2)$ for $\A_{even}$
gives a canonical embedding of
$\Proj(\A)$ (canonically isomorphic to $\Proj(\A_{even})$) into the
$\PP^2-$bundle $\PP(V_2)$.

Thus $\CCC := \Proj (\A)$
is a conic bundle  and we define $\pi_{\A} : \CCC \ra B$ as the restriction
of the natural projection  $\pi_2:\PP(V_2)
\ra B$.

Then the natural morphism
$\varphi_{\A}:S
\rightarrow \CCC = \Proj (\A)$ induced by the inclusion $\A \subset
\R$ yields a
factorization
$$f =\pi_{\A} \circ \varphi_{\A}.$$

Since the multiplication map from $V_3^+\otimes V_3^+$ to $\R_6$ has
image contained in $\A_6$ by lemma \ref{splittingR}, the ring
structure on $\R$ induces a map
$$ \delta:(V_3^+)^2 \rightarrow \A_6.
$$

\begin{df}\label{defP}
Let $P$ be a  point in the support of $\tau$.  We have seen in lemma
\ref{torsgenus2} that the map $\sigma_{2,P}$ has rank $2$, therefore
its image gives a pencil of lines in the plane which is the fibre of
$\PP(V_2)$ over this point.

This pencil of lines has a base point. Taking all the points thus associated
to the points of supp($\tau$) we get a subset of  $\PP(V_2)$ that
we will denote by $\P$. Note that the projection onto $B$ maps
$\P$ bijectively onto supp($\tau$).
\end{df}
\begin{teo}\label{ProjA}
$ \CCC = \Proj(\A) \subset \PP(V_2)$ is the divisor in the linear
system $$|\hol_{\PP(V_2)}(2) \otimes \pi_2^*(\det V_1)^{-2}| $$
with equation  induced
by the map $i_2$ defined in lemma \ref{structA}.

$\CCC$ has at most rational double points as singularities.

$\varphi_{\A}$ is the minimal resolution of the singularities of the
relative canonical model $ X : = \Proj(\R)$, a double cover of $\CCC$
whose branch
locus consists of the set of isolated points  $\P$ together with the divisor
$\Delta_{\A}$
   in the linear system
$|\hol_{\CCC} (3) \otimes \pi_{\A}^* (V_3^+)^{-2}
|$ determined by $\delta$ ($\Delta_{\A}$ is thus disjoint from $\P$).
\end{teo}
\Proof As observed in remark \ref{sigma2determinesA}, by the exact
sequence (\ref{exseqAeven}) $\A_{even}$ is the quotient of the
symmetric algebra $Sym(V_2)$ by the sheaf of principal ideals image
of the map $i_2:(\det V_1)^2 \rightarrow S^2(V_2)$, therefore its
$\Proj$ is a Cartier divisor $\CCC$ in the corresponding linear
system.

The map $\varphi_{\A}$ factors through $X = \Proj (\R)$. The natural map
$S \rightarrow X$ is the contraction of the $(-2)-$curves
   contained in a fibre of $f$; the map $ X \ra
\CCC$ is  finite of degree $2$ by lemma \ref{splittingR}; in
particular, because $X$ has only rational double points as
singularities, also the singularities of $\CCC$ must be isolated
and rational; since $\CCC$ is a divisor of the smooth $3-$fold
$\PP(V_2)$ they are rational hypersurface singularities, i.e.,
Rational Double Points (cf. \cite{artin}).

It only remains to compute the branch locus of the double cover $\psi: X
\ra \CCC$.

Since the question is local, we may assume that $ X $ is the
subscheme of
    $\PP(1,1,2,3) \times  B$ defined by the equations
$$ q (x_0, x_1, y, t) = 0, \  z^2 = g_6 (x_0, x_1, y, t).$$
Observe that $x_0 = x_1 =  y = 0$ implies then $ z=0$.

At a point where $x_i \neq 0$ we simply localize the two equations
dividing them by
$x_i^2$, respectively by $x_i^6 $. Hence, $z=0$ is the ramification
divisor and
$g_6 = 0$ is the branch locus. At the points where $ x_0 = x_1 = 0$,
$y=1$ we must have
a point of $\P$. This point is fixed for the involution $ z \ra -z$  and
$ g_6 \equiv  y^3 + x_0 \phi_5 + x_1 \psi_5 (mod \ t) $ does not vanish.
\qed

Theorem \ref{ProjA} shows that a
genus $2$ fibration is determined by a sheaf of algebras $\A$ constructed as in
lemma \ref{structA}, a line bundle $V_3^+$ and a map from
its square to $\A_6$. These data are first of all related to the invariants
of $S$
because of the exact sequence
$$ 0 \rightarrow Sym^2(V_1) \rightarrow V_2 \rightarrow \T_2
\rightarrow 0$$
which yields (see \cite{reid}) the Noether-type inequality
$$ (Horikawa)  \   K_S^2 -  6 (b-1) + 2 \chi = \deg (\tau) \geq 0.$$

   The following proposition shows that these data are not
independent even if we do not put restrictions on the invariants of $S$.
\begin{prop}\label{v3+}
$$V_3^+\cong\det V_1 \otimes \hol_B(\tau).$$
\end{prop}
\Proof We have decomposed the fibration $f$ as
$$ S \stackrel{r}{\rightarrow} X \stackrel{\psi}{\rightarrow}
\CCC \stackrel{\pi_{\A}}{\rightarrow} B,
$$ where $r$ is a contraction of $(-2)$-curves
to Rational Double Points, and $\psi$ is a finite double cover
corresponding to  the integral quadratic extension
$$ \A \subset \R \cong \A \oplus (\A[-3] \otimes V_3^+).$$

Localizing, we obtain the splitting
$$ \psi_* \hol_X \cong \hol_{\CCC} \oplus (\hol_{\CCC}(-3) \otimes
\pi_{\A}^*(V_3^+)).$$

Since $X$ and $\CCC$ have only Rational Double Points as singularities,
the standard adjunction formula for finite double covers yields
$$ \omega_X \cong \psi^* (\omega_{\CCC}(3 H) \otimes \pi_{\A}^*(V_3^+)^{-1})$$
where $H$ is a Weil divisor on $\CCC$ corresponding to $\hol_{\CCC} (1)$.
\
Therefore
$$ \psi^*  ( \hol_{\CCC}(1))\cong \psi^* (\omega_{\CCC |B}(3 H) \otimes
\pi_{\A}^*(V_3^+)^{-1})$$
and since by Theorem \ref{ProjA}
\begin{multline*}
\omega_{\CCC|B} =\pi_{\A}^* \det V_2 \otimes \hol_{\CCC}(-6)
\otimes \pi_{\A}^*(\det V_1)^{-2} \otimes
\hol_{\CCC} (4) =\\ = \hol_{\CCC}(-2)
\otimes \pi_{\A}^* (\det V_1 \otimes
\hol_B(\tau)),
\end{multline*}
we obtain the desired isomorphism by applying first $(\psi_*)^+$ and
then $(\pi_{\A})_*$.
\qed

Before stating the main theorem of this section we need to give an appropriate
definition.
\begin{df}\label{ass5ple}
Given  a genus $2$ fibration $ f : S \ra B$ we
define its {\bf associated 5-tuple} $(B,V_1,\tau,\xi,w)$
   as follows:
\begin{itemize}
\item $B$ is the base curve;
\item $V_1=f_*(\omega_{S|B})$;
\item $\tau$ is the effective divisor of $B$ whose structure sheaf is
     isomorphic to $\T_2$;
\item $\xi \in
       Ext^1_{\hol_B}(\hol_{\tau},S^2(V_1))/Aut_{\hol_B}(\hol_{\tau})$
       is the isomorphism class of the pair $V_2, \sigma_2$
$$ 0 \rightarrow S^2(V_1) \stackrel{\sigma_2}{\rightarrow} V_2
\rightarrow {\hol}_{\tau} \rightarrow 0;
$$
\item
setting
$$\tilde{\A}_6 : = \A_6 \otimes (V_3^+)^{-2},$$
$w \in \PP(H^0 (B, {\tilde{\A}_6 })) \cong |\hol_{\CCC}(6) \otimes
\pi_{\A}^*(V_3^+)^{-2}|$ is the class of a section with associated
divisor
   $\Delta_{\A}$.
\end{itemize}
\end{df}
\begin{df}\label{F}
Let $B$ be a smooth curve, $V_1$, $V_2$ two vector bundles on $B$ of
respective ranks $2$ and $3$. Let $\sigma_2:S^2(V_1) \rightarrow V_2$
be an injective homomorphism whose cokernel is isomorphic to
the structure sheaf of an effective
divisor $\tau$ on $B$.

We define ${\tilde{\A}_6 }$ to be the vector bundle
$$(coker\ i_3) \otimes \left( \det V_1 \otimes \hol_B(\tau)
     \right)^{-2},$$ where the map $i_3:(\det V_1)^2 \otimes V_2
     \rightarrow S^3(V_2)$ is the one induced by $\sigma_2$ as in lemma
     \ref{structA} (see exact sequence (\ref{exseqAeven})).
\end{df}
\begin{oss}
The fact that ${\tilde{\A}_6 }$ is automatically a vector bundle follows by the
assumption that the rank of $\sigma_2$ drops at most by $1$, which
implies that the rank of the induced map $i_3$ never drops: this is a
straightforward computation that we leave to the reader.
\end{oss}
The previous results and definitions allow us to introduce the
building package of
a genus $2$ fibration:
\begin{df}\label{adm5ple}
We shall say that a a $5-$tuple
    $(B,V_1,\tau,\xi,w)$ is an {\bf admissible genus two $5$-tuple} if
\begin{itemize}
\item $B$  a smooth curve;
\item $V_1$  a vector bundle on $B$ of rank $2$;
\item $\tau \in Div^+(B)$;
\item $\xi \in
Ext^1_{\hol_B}(\hol_{\tau},S^2(V_1))/Aut_{\hol_B}(\hol_{\tau})$
yields a vector bundle $V_2$;
\item $w \in \PP(H^0 \left(B, {\tilde{\A}_6 } \right))$, where
${\tilde{\A}_6 }$ is
the vector bundle determined by $\xi$ as in definition \ref{F},
\end{itemize}
and moreover the following (open) conditions are satisfied:
\begin{itemize}
\item[i)] Let $\A$ be the sheaf of algebras determined by
     $B,V_1,\tau,\xi$ as in remark \ref{sigma2determinesA}; then the 
conic bundle
     $\CCC : = \Proj(\A)$ has at most rational double points as singularities.
\item[ii)]  Let $\Delta_{\A}$ be the divisor of $w$ on
     $\CCC$; then $\Delta_{\A}$ does not contain
    any point of the set $\P$ defined as in  \ref{defP}.
\item[iii)] At each point of
     $\Delta_{\A}$ the germ  of the double cover $X$ of
     $\CCC$ branched on $\Delta_{\A}$  has at most rational
     double points as singularities.
\end{itemize}
\end{df}
\begin{teo}\label{gen2}
  Let $f$ be a relatively minimal genus $2$ fibration. Then its
{\bf associated 5-tuple} $(B,V_1,\tau,\xi,w)$ is admissible.

Viceversa, every admissible genus two $5-$tuple determines a
sheaf of algebras $ \R \cong \A \oplus ( \A [-3] \otimes \det(V_1) 
\otimes \hol_B
(\tau))$ over $B$
whose relative projective spectrum $X$ is the relative canonical model
of a relatively minimal genus $2$ fibration $f:S \rightarrow B$
having the above as associated 5-tuple. Moreover, the
surface $S$ has
the following invariants:
$$\begin{array}{lll} \chi(\hol_S)&=&\deg(V_1)+(b-1)\\ K_S^2&=&2\deg V_1
+ \deg \tau + 8(b-1).
\end{array}$$

We obtain thus a bijection between (isomorphism classes of)
  relatively minimal genus $2$
fibrations and  (isomorphism classes of) associated
$5-$tuples are isomorphic, which is functorial in the sense that to a flat family of
fibrations corresponds a flat family of 5-tuples.
\end{teo}
\Proof  The associated $5-$tuple of a genus $2$ fibration
is admissible by theorem \ref{ProjA}.

Viceversa, given an admissible (genus two) $5-$tuple, let
$\sigma_2:S^2(V_1) \rightarrow V_2$ be the map induced by $\xi$ and
  construct a sheaf of algebras $\A$ as in remark
\ref{sigma2determinesA}.

Finally, a representative of $w$ provides the sheaf of modules $\A \oplus
\left( \A[-3] \otimes \det V_1 \otimes \hol_B(\tau) \right)$
with the structure of a sheaf of algebras $\R$.

We get then by construction a finite morphism  of degree $2$
$$\psi: X := \Proj(\R)
\rightarrow
\CCC : = \Proj(\A)$$  whose branch locus is the union of
the divisor $\Delta_{\A}$ and of the finite set $\P$, which are
disjoint by condition ii), while condition iii) ensures that $X$ has at most
rational double points as singularities. We finally take $S$
to be a minimal resolution of the
singularities of $X$.

The induced map $f: X \rightarrow B$ is clearly a genus $2$
fibration whose associated $5-$tuple coincides with the given
admissible genus two $5-$tuple:
    the only non trivial verification,  that $f_*
\omega_{S|B}=V_1$, can be done by the same arguments
used in proposition \ref{v3+}. $f$ is relatively
minimal since by construction the relative canonical
bundle is  $f-$nef.

If we start from  relatively minimal genus $2$ fibration $f:S \rightarrow B$
we get it back from the associated 5-tuple, and the invariants
are as stated because of \ref{fujitadeg}.

Finally, given a flat family $ f_T : \SSS \ra  \B  \ra T$ of genus $g$ fibrations,
the sheaves $ \V_{T,n} : = (f_T )_* (\omega_{\SSS |  \B }^{\otimes n})$ are
vector bundles on $\B$ because of the base change theorem (Thm. 12.11
of \cite{har}). For this reason the exact sequence
$$ 0 \ra S^2(\V_{T,1}) \stackrel{\sigma_{T,2}}{\longrightarrow }\V_{T,2} \ra \T_{T,2}
\ra 0$$ remains exact when restricted to a fibre $ B_t$, and 
since the rank of $\sigma_{T,2}$ drops at most by one, $\T_{T,2}$
is the structure sheaf of a Cartier divisor on $\B$. The rest follows
now easily. 
\qed

\begin{oss}\label{condareeasy}
We discuss here the open conditions in definition \ref{adm5ple}.

The fibres of the conic bundle $\CCC$ over a point not in
supp$(\tau)$ are smooth conics, whereas the fibres over points of supp
$(\tau)$ are singular: more precisely a singular reduced conic if the
kernel of the map $\sigma_2$ on the corresponding fibre is not a
square tensor,  a double line otherwise.

Hence $\CCC$ is  smooth on the complement of the special fibres and
   on the nonsingular locus of $\CCC$
the only condition to be fulfilled is that
$\Delta_{\A}$ have only simple singularities (cf. \cite{bpv}).

On a reduced singular fibre there is only one singular point, namely,
the point $ P(b) \in
\P$, and a local computation shows that $P(b)$ is automatically a 
singularity of type
$A_{2s+1}$: we must further  check that
$ P(b) \notin \Delta_{\A}$.

The situation is slightly more complicated when $\CCC$ has a double
fibre $F_b$ : one has first to check that
$F_b$ is not contained in the singular locus. If this is the case,
writing the equation locally one finds two possible cases:

case i) $\CCC$  has only a singular point (of type
$D_{2s}$) in $P(b)$, which, as before,
must not lie in $\Delta_{\A}$.

case ii) $P(b)$  is a singularity of type $A_1$, and $\CCC$
has on the fibre $F_b$ a further singular point $P'$, still of type $A_1$.
The branch curve $\Delta_{\A}$ is allowed to pass through $P'$,
  and if this is the case one has to verify that the double cover $X$
has a
Rational Double Point.

We  have thus also seen that the only singularities
that $\CCC$ can have are of type $A_{2s+1}$ and of type $D_{2s}$.
\end{oss}

\begin{oss}
In \cite{horp} Horikawa gave a classification of the special fibres of
a genus $2$ fibrations. Using the local equations we have used
throughout all this section it is not difficult to recognize the
classification of Horikawa, and more precisely how the parameters $s$
and $\lambda$ we have introduced determine the geometry of the
fibre. We do not enter in the detail of the straightforward but long
computation, but we think it can be of some interest for the reader to
know that the special fibre is (in the notation of \cite{horp}) of
type:
\begin{eqnarray}
\nonumber I_{\frac{s+1}{2}} & {\rm if }\ s \ {\rm is\ odd \ and }&
\lambda\neq 0;\\ \nonumber II_{\frac{s}{2}} & {\rm if }\ s \ {\rm is\
even \ and }& \lambda \neq 0;\\
\label{tabella} III_{1}\ {\rm or }\ V & {\rm if }\ s =1 \ {\rm and }&
\lambda=0;\\
\nonumber III_{\frac{s+1}{2}} & {\rm if }\ s \geq 3 \ {\rm is\ odd\
and }& \lambda=0;\\ \nonumber IV_{\frac{s}{2}} & {\rm if }\ s \ {\rm
is\ even \ and }& \lambda=0.
\end{eqnarray}
\end{oss}
\section{The genus of the bicanonical pencil of a Godeaux surface}
We give an application of theorem \ref{gen2} to the classification of
Godeaux surfaces. The following theorem was already proved by Bombieri
(\cite{bomman}, cf. also footnote (1) on page 494 of \cite{bom})
using Ogg's list of genus $2$ fibres (cf.
\cite{ogg}).
\begin{teo}\label{godeaux}
Let $S'$ be a minimal numerical Godeaux surface, i.e., the minimal
model of a surface
with
$ K_{S'}^2 = 1$, $p_g (S') =0$, and let
$f : S \rightarrow
\PP^1$ be the fibration induced by the bicanonical pencil of $S'$. 
Then the genus
of the fibre can only be $3$ or $4$.
\end{teo}
\Proof
We already know (cf. \cite{cp} and \ref{minimalgen2}) that it
suffices to exclude the
case where $S=S'$ and $|2K_S|=| M| + \Phi $ where $|M|$ is a base point free
genus $2$
pencil ($K_S\Phi=M^2= 0, K_S M = 2$).

Let's argue by contradiction: we have $b=0$, $K_S^2=\chi(\hol_S)=1$.
By the formulae in
theorem
\ref{gen2} we conclude $\deg V_1=2$, $\deg \tau=5$. By proposition
\ref{v3+} we get $V_3^+=\hol_{\PP^1}(7)$.

In particular
\begin{multline*}
h^0(\hol_{S}(K_S+\Phi))=h^0(\hol_{S}(3K_S-M))=h^0(f_*\hol_{S}(3K_S-M))=\\
=h^0(f_*(\omega_{S|B}^3)\otimes \hol_{\PP^1}(-7))\geq h^0(V_3^+\otimes
\hol_{\PP^1}(-7))=h^0(\hol_{\PP^1})=1.
\end{multline*}

By the long cohomology exact sequence associated to
$$ 0 \rightarrow \hol_S(K_S) \rightarrow \hol_S(K_S+\Phi) \rightarrow
\hol_\Phi(K_S+\Phi) \rightarrow 0,
$$ and since $p_g(S)=q(S)=0$, we conclude that $h^0(\hol_\Phi(K_S+\Phi))\geq
1$.

On the other hand  the long cohomology exact sequence associated to
$$ 0 \rightarrow \hol_S \rightarrow \hol_S(\Phi) \rightarrow
\hol_\Phi(\Phi) \rightarrow 0,
$$ implies $h^0(\hol_\Phi(\Phi))=0$ since  $q(S)=0$.

Since $K_S \Phi=0$ we know  (cf. \cite{artin}) that
we have an isomorphism
$\hol_\Phi(K_S)  \cong \hol_\Phi$, and we derive a contradiction from
$h^0(\hol_\Phi(K_S+\Phi))\geq 1$, $h^0(\hol_\Phi(\Phi))=0$.
\qed
Observe, by the way, that the above proof does not use the full
strength of
theorem \ref{gen2}.
\section{Surfaces with $p_g=q=1$ and $K^2=3$}
Let $S$ be a minimal surface of general type with $p_g=q=1$. 
The Noether inequality
and  the Bogomolov Miyaoka Yau inequality
give us
$2 \leq K^2_S \leq 9$. The Albanese map is a morphism $f:S
\rightarrow B$ where
$B$ is a smooth elliptic curve.

Theorem \ref{gen2} can be used to study the case in which the genus $g$
of the general fibre of $f$ is $2$. Xiao's inequality (\cite{xiao}, thm 2.2)
in this case yields $2 \leq K^2_S \leq 6$. In fact,\cite{k^2=2} proved that
 for $K^2_S=2$ the genus of the Albanese fibre is $2$, the
surface is a double cover of  the second symmetric product $B^{(2)}$ of $B$,
   and the moduli space is irreducible, generically smooth, unirational of
dimension $7$.

The class of surfaces of general type with $K^2=3$, $p_g=q=1$ was
investigated in \cite{cc1}, \cite{cc2}. \cite{cc1}  proved that for
this class of surfaces $g=2$ or $3$, and the second case is completely
classified in \cite{cc2}, where it is shown that the
corresponding component of the moduli space is
generically smooth, irreducible, unirational of dimension $5$.

\cite{cc1}  shows that the surfaces with $p_g=q=1$, $K^2=3$ and genus
$2$ of the Albanese fibre are double covers of $B^{(2)}$.
Given $ t \in B$,  let $D_t$ be the curve of divisors
in $B^{(2)}$ containing $t$, i.e.,  $D_t : = \{ ( x + t) |  x \in B \} $,
and let $E_t$ be the fibre over $t$ ($\cong \PP^1$) of the Abel Jacobi map
$B^{(2)} \ra B.$
In order to show that this  class
of  surfaces is not empty, it was shown in loc. cit.  the existence of curves
$C  \sim 6 D - E$  ($C$ algebraically equivalent to $6D-E$)  and with exactly two ordinary
triple points on a given fibre $E$. By taking as branch locus the union $ C
\cup E$, and a corresponding double cover, one obtains  surfaces with
$p_g=q=1$, $K^2=3$ and genus
$2$ of the Albanese fibre.
   It was conjectured there
(see problem 5.5) that this family of surfaces should form an irreducible
family of the moduli space: we will now disprove  this conjecture
   using theorem \ref{gen2} , and showing moreover that the corresponding
locus of the moduli space is
indeed disconnected.

\begin{teo}\label{Kquadro3gen2}
Let  $\M$ be the moduli space of the minimal surfaces
of general type  $S$ with $p_g(S)=q(S)=1$,
$K_S^2=3$. Then $\M$ has $4$ connected components, all irreducible and unirational
of dimension $5$.
\end{teo}
\begin{oss}
The unirational family of dimension $5$ studied in \cite{cc2} is  a connected
component of $\M$ by virtue of
the differentiable invariance of the genus $g$ (cf. \cite{isogenous}).
\end{oss}
Let $\M'$ be the (open and closed) subset of $\M$ corresponding
to the surfaces whose Albanese fibres have   genus $2$;
 theorem \ref{Kquadro3gen2} follows then from the
following:
\begin{prop}
   $\M'$ has $3$ connected components, all irreducible
and unirational of dimension $5$.
\end{prop}
To prove this we define a stratification of $\M$;  we first describe
in geometric terms one of the
strata, that we call {\bf of the main stream}, because
it is defined by generality assumptions on the branch curve.

The conjecture of \cite{cc1} was in fact related to the existence and
density of this
main stream family; however,   the surfaces constructed in \cite{cc1} were
not of the main stream!

The main stream family corresponds to double covers of $B^{(2)}$ branched on
curves $C + E$  belonging to the $5$-dimensional family 
(here $B$ also varies) described as follows
\begin{teo}\label{mainstream}
Let $ p : B^{(2)} \ra  B$ be the Albanese (Abel-Jacobi) map, and let $E$ be
a fixed fibre of $p$. Let $D_t$ be a natural section of $p$ given by the set
of pairs (of points) on $B$  which contain a given point $t \in B$.

Let $P_1, P_2$ be general points of $E$ and let $C$ be a divisor
algebraically equivalent to $ 6 D - E$ having $P_1, P_2$  as points of
multiplicity at least $3$.

Then, for $C$ general in the algebraic equivalence class, the linear system
$|C|$ has dimension $1$, the general curve $C$ inside the system is
irreducible and its proper transform after the blow up of $P_1$ and $P_2$
has singularities at most double points.
\end{teo}
\Proof
Recall (cf. e.g.  \cite{k^2=2}, \cite{cc1}, \cite{cc2}) that on  $B^{(2)}$ the
bianticanonical system
$|-2K|$ is a base point free pencil of elliptic curves $\Delta_s = \{
\{x,x+s\}\}$ where $s\in B$.
$\Delta_s $ is an irreducible smooth elliptic curve isomorphic to $B$, except
when $s$ is a non trivial
$2$-torsion point $\eta$, and then $\Delta_{\eta} = 2 T_{\eta}$, where
$T_{\eta}$ is the smooth elliptic curve $\cong B/ \langle\eta \rangle$.
Observe moreover that $-K$ is algebraically equivalent to $2 D - E$.

CLAIM 1)  Let $C$ be algebraically equivalent to $ 6 D -
E$: then   $$H^1(\hol_{B^{(2)}}(C)) = 0.$$
   {\em Proof of 1)}: by Serre duality it suffices to show the vanishing
   of
$H^1(\hol_{B^{(2)}}(K-C)) $.

Writing $H^1(\hol_{B^{(2)}}(K-C)) =
H^1(\hol_{B^{(2)}}(-(C-K)))$ we infer the vanishing from the Ramanujam
vanishing theorem ( \cite{ram}, Theorem 2, page 48, see also \cite{bom}). In fact, we
have the following algebraic
equivalences:
$ C - K \sim (-2K) + 4D \sim \Delta_s + D_{t_1} \dots + D_{t_4} $ which show
that indeed $C-K$ is linearly equivalent to a reduced and connected
effective  divisor which is not composed of an irrational pencil.
\QED for 1)

By Riemann Roch we obtain $\dim |C| = \frac{1}{2} C (C-K) - 1 = 13 $, whence
   for any choice of distinct points $P_1, P_2 \in B^{(2)}$ we have:

2) $\dim \ | C - 3 P_1 - 3 P_2 | \geq 1.$

3) Let $P_3$ be a third point of the fibre $E$, distinct from $P_1, P_2$;
 a curve  in $| C - 3 P_1 - 3 P_2 |$ passing through $P_3$ must
contain the whole fibre $E$, whence the system $| C - 3 P_1 - 3 P_2 |$
contains  $ E +| C - E - 2 P_1 - 2 P_2 |$ as a linear subsystem  of
codimension at most 1.

4) Observe that  $ \dim |C - E| = 6 $ follows by the  same argument as in 1).
Instead, we claim that $ C - 2 E \sim - 3 K$ is empty if the linear
equivalence class of $C$ is general.
In fact, if $M \sim (C-2E)$ is an effective divisor, then $ M \Delta_s = 0$,
thus $M$ consists of a fibre $\Delta_s$ of the bianticanonical pencil plus a
curve $ T_{\eta}$.

5) Assume now that the general curve in  $ | C - 3 P_1 - 3 P_2 |$ is
always reducible, i.e., that generally  there exists   a fixed part
$\Phi$.

There are two cases:

5.1) $\Phi \not\geq E$.

5.2) $\Phi \geq E$.

CLAIM I): {\em case 5.1) is not possible.}

{\em Proof of claim I)}: If  $\Phi \not\geq E$, by 3) 
the system $ | C - 3 P_1 - 3 P_2 |$
would contain
a reducible curve of the form $ E + \Phi + M$, where $C' : = \Phi + M$ is in
$ | C' - 2 P_1 - 2 P_2 |$.  Now let the general pair specialize to a pair
of distinct points contained in the complete intersection
$ \Delta_s \cap E$.

The intersection number $ C' \Delta_s $ equals $ \Delta_s \cdot 2 D = 4$,
and therefore
the  system  $ | C' - 2 P_1 - 2 P_2 |$ contains the subsystem $\Delta_s +
| C' - \Delta_s -  P_1 -  P_2 |$ as a linear subsystem of codimension
$\leq 1$.

Since $C' - \Delta_s \sim 2 D$, and $ |2 D - E| = \emptyset$ for general choice
of the linear equivalence classes of $D,E$, we
see that $ \dim |2D|= 2$ always, and that for general choice of
the linear equivalence class of $C'$, the subsystem
$| C' - \Delta_s -  P_1 -  P_2 |$  has dimension zero.

MAIN CLAIM: {\em For special choice of $\{P_1, P_2 \} \subset \Delta_s \cap E$
and  $C$, $\Delta_s$ general,
$ | C' - 2 P_1 - 2 P_2 |$ equals the subsystem $\Delta_s +
| C' - \Delta_s -  P_1 -  P_2 |$. Moreover, for a general choice of
   $\{P_1, P_2 \} $ and of $C'$, $ | C' - 2 P_1 - 2 P_2 |$ has
   dimension $0$.}

{ \em Proof of the Main Claim.}
For otherwise,  $|C' - 2 P_1 - 2 P_2 |$  always has dimension
$\geq 1$ for a special choice of the points $\{P_1, P_2 \} \subset
\Delta_s \cap E$,
and by 4) this dimension is  generally equal to $1$ .
Varying $E$, $\Delta_s$, and the linear equivalence class, we obtain an
irreducible family of dimension $4$ which dominates a subvariety
of dimension $2$ of the symmetric square of $B^{(2)}$.

Thus, it suffices to show that at some curve $Z$ of the $4$-dimensional
family  given by curves in $\Delta_s + G$ with $  G \in
| C' - \Delta_s -  P_1 -  P_2 |$, the equisingular deformations of $Z$
dominate the symmetric square of $B^{(2)}$ and the dimension of the
fibre is at most $1$. In fact we then  obtain a contradiction which
proves the first assertion. Since the family has dimension  $ \geq 5$,
the second assertion also follows.

Let us choose $Z$ of the form $\Delta_s  + D_1 + D_2$, where $P_i \in
D_i$.

Observe that $D_1 $ and $D_2$ meet transversally at a point $P$, while
$\Delta_s  \cap D_i = \{ P_i, Q_i\}$.

We make everything explicit, writing 
$ E = \{ (y, -y) | y \in B\}$ for the fibre of pairs that add to
$0$, and
 letting $P_1 = (a,-a), \ P_2 = (b,-b)$.
Set $D_1 = \{ (a,x)| x \in B\}$, $D_2 = \{ (b,z)| z \in B\}$, $s = 2a$.
Then $\Delta_s = \{ (u , u + 2a)\}$. Without loss of generality we may assume
that $ 2a = s= 2b$. Then one computes that $Q_1 = (a, 3a)$, while
$Q_2 = (b, 3b)$.

We conclude that on the elliptic curve $\Delta_s$ the divisor $P_1 + P_2 $
is not linearly equivalent to $ Q_1 + Q_2 $ for a general choice of $s$, since
$ -a -b \neq a + b$.

Consider now the curve $Z$ and the normal sheaf $ \NNN'_Z $ of
deformations that are equisingular  at  $P_1 , P_2 $.

Its space of global sections fits into the exact sequence:
\begin{multline*} 0 \ra H^0(\NNN'_Z ) \ra \\
\ra H^0(\NNN_{D_1 }( Q_1 + P)) \oplus
H^0(\NNN_{D_2}( Q_2 +  P)) \oplus H^0( \hol_{\Delta}(Q_1 + Q_2)) \ra\\
\ra \C_{Q_1}\oplus \C_{Q_2}\oplus \C_{P} \ra
\end{multline*}

Since the normal sheaf  $\NNN_{D_i}$ has degree $1$,
each summand  $H^0(\NNN_{D_i }( Q_i + P))$ in the middle
surjects onto the direct sum $\C_{Q_i}\oplus \C_{P}$ of the  two
corresponding  summands on the right.
It follows  that  we have surjectivity at the right term, and hence
that $ \dim \ (H^0(\NNN'_Z ) ) =  5$.

Since we saw that
$P_1 + P_2 $
is not linearly equivalent to $ Q_1 + Q_2 $ on $\Delta_s$, we obtain that
$ H^0(\NNN'_Z )  $ surjects onto the direct sum of the cotangent spaces of
the surface at $P_1, P_2$, since $ H^0(\NNN_{D_1 }( Q_1 + P -P_1)) \oplus
H^0(\NNN_{D_2}( Q_2 +  P - P_2)) \oplus H^0( \hol_{\Delta}(Q_1 + Q_2 - P_1 -
P_2)) $ has dimension $4$ and surjects onto $ \C_{Q_1}\oplus
\C_{Q_2}\oplus \C_{P}$.

We have therefore shown that the tangent dimension to the family is
$5$ at $Z$, the
differential of the map onto the symmetric square of $ B^{(2)}$ is surjective,
whence the family is smooth at $Z$ of dimension $5$ and the fibre is
smooth of dimension $1$.
\QED for the Main Claim.

We now proceed  to prove Claim 1): we let the pair $P_1, P_2$ tend
to a special pair: then the reducible curve $\Phi + M$ has as limit a curve
in the linear system $| C' - 2 P_1 - 2 P_2| = \Delta_s + G $.
By a general choice of the class of $C'$ we obtain that $G \sim 2 D$ is
irreducible. It follows that, modulo exchanging the roles of $\Phi$ and
$M$, we may assume that  $\Phi$ tends to $\Delta_s$, in particular
$\Phi$ is algebraically equivalent to $\Delta_s$ and thus belongs to the
antibicanonical pencil.
This however shows that the points $P_1, P_2$ must always be
special, a contradiction.
\QED for claim 1)

We finally observe that also case 5.2) is impossible: in this case we
would have $|C- 3 P_1 - 3 P_2| = E+ |C- E- 2 P_1 - 2 P_2|$ and then
 $\dim |C- E- 2 P_1 - 2 P_2|$ is always $\geq 1$, contradicting the
Main Claim.

Hence, we have proved that the general curve $C$ which has multiplicity at
least $3$ in each point $P_1, P_2$ is generally irreducible, and for a fixed
choice of $P_1, P_2$ we have a linear pencil.

Observe however that the intersection number of $C$ with $E$ is $6$, whence
it follows that each point $P_i$ is of multiplicity exactly $3$, and that the
fibre is not tangent. After blowing up the two points, the proper transform
of $C$, call it $ \tilde{C}$, does not intersect the proper transform of the
fibre $E$, and moreover  $ \tilde{C}^2 = 6$, whence we obtain a pencil of
curves  $ \tilde{C}$ which generically have at most double points as
singularities.
\QED

We have just shown how to construct the main
stream family of surfaces, i.e., the double covers for which all the
choices made (elliptic curve, fibre $E$, two points on the fibre,
linear equivalence class) can be taken general. However, for special
choices, it can happen that the dimension of the corresponding linear
system of curves jumps. 

It is for instance clear that this happens when the points are very special,
i.e., on a curve $T_{\eta}$, but  to analyse all the possible cases,
we  need to make full use of the theory developed in the previous
section.

We need now to recall some standard notation and results about vector bundles
over elliptic curves.
\begin{df}
Given a point $u$ of an elliptic curve $B$, and integers $r,d$ with $r>0$,
$(r,d)=1$, we will denote by $E_u(r,d)$ (following \cite{atiyah}) the only
indecomposable vector bundle of rank $r$ on $B$ with $\det E_u(r,d)
=\hol_B(u)^{\otimes d}$.
\end{df}
\begin{lem}[\cite{cc1}, page 76]\label{cirofab}
Let $S$ be a minimal surface of general type with  $p_g(S)=q(S)=1$
whose Albanese map $f:S \rightarrow B$ is a genus $2$ fibration.

If $K^2_S \leq 3$, then there is a point $u \in B$ such that
$f_*\omega_{S|B}=E_u(2,1)$.

If instead $K^2_S \geq 4$, either $f_*\omega_{S|B}=E_u(2,1)$ or
$f_*\omega_{S|B}=L
\oplus \hol_B(u)$ with $u \in B$ and $L$ a non trivial torsion line bundle.
\end{lem}

From this lemma and theorem \ref{gen2} we immediately get the following
\begin{oss}\label{trasloin0}
Let $S$ be a minimal surface of general type with  $p_g(S)=q(S)=1$ and
$K^2_S=3$ whose Albanese map $f:S \rightarrow B$ is a genus $2$ fibration.
Let $(B,V_1,\tau,\xi,w)$ be the associated $5-$tuple:
then we may assume w.l.o.g. that
$V_1=E_{[0]}(2,1)$.

Moreover $\deg \tau= K^2_S-2=1$, i.e., $\tau$ is a point of  $B$.
\end{oss}

To apply theorem \ref{gen2}, we need to compute
$S^2(V_1)=S^2(E_{[0]}(2,1))$; this is a computation based on the results
in \cite{atiyah}: the interested reader can find the computation in
\cite{cc2}.
\begin{oss}\label{S2(E)}[cf. \cite{cc2}]
Given an elliptic curve $B$, let us denote by $L_i,\ i\in \{1,2, 3\}$ the
three line bundles on $B$ with $L_i \not\cong \hol_B$, $L_i^2 \cong
\hol_B$. Then
$$S^2(E_{[0]}(2,1))=\bigoplus_{i=1}^3 L_i([0]).$$
\end{oss}
The $4$th element $\xi$ of the associated $5-$tuple belongs to
$$Ext^1_{\hol_B}(\hol_{\tau},
\oplus_{i=1}^3 L_i([0])
)/Aut_{\hol_B}(\hol_{\tau}).$$

Note that $\xi$ can't be the class of $0$, since the extension must
yield a vector bundle. Since $\tau$ is a point,
$Aut_{\hol_B}(\hol_\tau) \cong \C^*$, acting on
$Ext^1_{\hol_B}(\hol_{\tau},S^2(V_1)) \setminus \{0\}$ by scalar
multiplication: therefore the space to which $\xi$ belongs is the
projective space $\PP(Ext^1_{\hol_B}(\hol_{\tau},S^2(V_1)))$.

To compute this space we fix a section  $f_{0}\in H^0(\hol_B(\tau))
   \setminus \{0\}$; applying the functor
   $Hom_{\hol_B}(\cdot,S^2(V_1))$ to the exact sequence
$$ 0 \rightarrow \hol_B([0]-\tau) \stackrel{\cdot (-f_{0})}{\longrightarrow}
\hol_B([0]) \rightarrow \hol_{\tau} \rightarrow 0,
$$
we get isomorphisms
\begin{multline}\label{calcolaxi}
Ext^1_{\hol_B}(\hol_{\tau},\oplus_{i=1}^3 L_i([0]))
\cong
Hom_{\hol_B}(\hol_B([0]-\tau),\oplus_{i=1}^3 L_i([0]))
\cong\\
\cong
H^0(\oplus_1^3 L_i(\tau))
\cong
\C^3
.
\end{multline}
Therefore the associated projective space is a $\PP^2$: we have
parametrized the first $4$ data of our $5-$tuples via a
unirational parameter space of dimension $4$ (a $\PP^2-$bundle over a
universal family of  elliptic curves).

To complete the $5-$tuple we need to compute $H^0(\tilde{\A}_6)$ (in terms
of $B,\tau,\xi$).  First we describe $V_2$ explicitly as follows
\begin{oss}
Let $(f_1,f_2,f_3)$ be an element of $H^0(\oplus_1^3 L_i(\tau))$. Then the
corresponding extension class $\xi$ induced by $f_0$ as above
 yields a sheaf
$V_2$ that is the cokernel of the map from $\hol_B([0]-\tau)$
to $\hol_B([0]) \oplus \left( \oplus_1^3 L_i([0]) \right)$ induced by
the matrix $ ^t(f_0,f_1,f_2,f_3)$.
\end{oss}
In fact, by standard homological algebra we have a commutative diagram
$$
\begin{array}{ccccccccc}
0
&
\rightarrow
&
\hol_B([0]-\tau)
&
\stackrel{\cdot (-f_{0})}{\rightarrow}
&
\hol_B([0])
&
\rightarrow
&
\hol_{\tau}
&
\rightarrow
&
0\\
&&\downarrow&&\downarrow&&||&&\\
0
&
\rightarrow
&
\oplus_1^3L_i([0])
&
\stackrel{\sigma_2}{\rightarrow}
&
V_2
&
\rightarrow
&
\hol_{\tau}
&
\rightarrow
&
0,
$$
\end{array}
$$
where the vertical map on the left is the one given by
${}^t(f_1,f_2,f_3)$. This diagram clearly induces an exact sequence
\begin{equation}\label{dafaV2}
0
\rightarrow
\hol_B([0]-\tau)
\rightarrow
\hol_B([0]) \oplus (\oplus_1^3L_i([0]))
\rightarrow
V_2
\rightarrow
0.
\end{equation}
\begin{lem} \label{classV2}
$V_2$ is determined by the $f_i$'s as follows:
$$
\begin{tabular}{llll}
I) & $ f_i \neq 0 \ \forall i\in \{1,2,3\} $ & $\Leftrightarrow$& $V_2(-[0])
\cong E_{\tau}(3,1)$\\
II) & $\exists !i\in \{1,2,3\}$ with $f_i = 0$& $\Leftrightarrow$& $V_2(-[0])
\cong E_{\tau_i}(2,1)\oplus L_i$\\
III) & $\exists !i\in \{1,2,3\}$ with $f_i \neq 0$& $\Leftrightarrow$&
$V_2(-[0]) \cong L_i(\tau) \oplus L_j \oplus L_k$\\
\end{tabular}
$$ where the point $\tau_i \in B$ is the divisor of a non trivial section of
$L_i(\tau)$, and whenever $j$ and $k$ appear $\{i,j,k\}=\{1,2,3\}$.\\
\end{lem}
\Proof
Let $m$ be the cardinality of the set $\{i| f_i=0 \}$. Note that
$m\leq 3$ because $f_0$ is different from zero by assumption.

By the exact sequence (\ref{dafaV2}), follows that $V_2$
is a vector bundle if and only if the
$f_i$'s have no common zeroes, i.e., if and only if $m\leq 2$
(since the points $\tau_i$ are distinct).

Tensoring the exact sequence (\ref{dafaV2})
by $\hol_B (-[0])$  and since
$$ H^1(\hol_B(-\tau) )\stackrel{ \cdot f_0}{\rightarrow}  H^1(\hol_B)$$
is an isomorphism,
we conclude that $H^1(V_2(-[0]))=0$.

Further twisting the exact sequence (\ref{dafaV2}) by any degree $0$ line
bundle $L$ and repeating the argument, we see that
$H^1(V_2(-[0])\otimes L)$ vanishes unless $L=L_i$ for some $i$ and $f_i=0$: in
this last case $H^1(V_2(-[0])\otimes L) \cong \C$.

Therefore  $V_2(-[0])$ is a vector bundle of rank $3$
and determinant $\hol_B(\tau)$ having the property that there are
exactly $m$ line bundles $L$ of degree zero with $H^1(V_2(-[0])\otimes
L) \neq 0$.

If $V_2(-[0])$ is indecomposable, it is $E_{\tau}(3,1)$ by Atiyah's
classification. This sheaf has trivial cohomology when twisted with
any degree $0$ line bundle, therefore $m=0$: we are in case I).

If on the contrary $V_2(-[0])$ is a sum of three line bundles, being
   a quotient of $\hol_B \oplus (\oplus L_i)$ of degree $1$, we see that
two summands have degree $0$: it is clear that in this case $m=2$ and
that $V_2(-[0])$ is the one described in the statement (case III)).

Else $V_2(-[0])=W \oplus L$ is the sum of two indecomposable vector bundles of
respective ranks $2$ and $1$.

First we exclude the case $\deg W=0$, $\deg L=1$.  Recall that, by
theorem $5$ of \cite{atiyah}, part ii), for each indecomposable vector
bundle $W$ of degree zero there is exactly one line bundle $L$ with
$H^1(W \otimes L) \neq 0$. Therefore, if $W$ had degree $0$, we would get
$m=1$. But then, if $f_i$ is the vanishing section, $L_i$ is a direct
summand of $V_2(-[0])$, a contradiction.

Then (by semipositivity) $\deg L=0$ and $\deg W=1$; every twist of $W$
by a degree $0$ line bundle has trivial first cohomology group,
therefore the corresponding twist of $V_2(-[0])$ has nontrivial first
cohomology group if and only if we twist by $L^{-1}$: whence $m=1$ and we are
   in case II).
\QED

The preceding result suggests to consider a stratification of our
   moduli space of surfaces:
\begin{df}\label{strata}
We stratify $\M'$ as $\M'= \M_I \cup \M_{II} \cup \M_{III}$
according to the number  of  indecomposable summands for
$V_2=f_*(\omega_{S|B}^2)$,
as in lemma \ref{classV2}.

We decompose  $\M_I$ further as
$\M_I^{\circ}
\cup
\M_{I,3}
$
where $\M_{I,3}$ consists of the surfaces with $3\tau \equiv 3[0]$
and $\M_I^{\circ}$ is the rest.
\end{df}
\begin{oss}\label{strategy}
By a theorem of Clemens (\cite{clem}, see  \cite{fabfrank} for the way
it is applied) the  dimension of an irreducible component of the
moduli space of minimal surfaces of general type with $p_g=q=1$ is
at least $10\chi-2K^2+1$, that is, $5$ in our case. In particular each
stratum of smaller dimension cannot contain an open set and can be disregarded
for the determination of the irreducible components.

Studying the   $4$ strata introduced above we will find that $\M_{III}$
consists of two
not empty unirational families of dimension $5$, $\M_I^{\circ}$ is
the unirational main stream
family, while
   the other strata have smaller dimension.

Finally we will study how the closures of these strata 
intersect.
\end{oss}
To simplify the exposition, we  first treat together the cases in
which $V_2$ is decomposable: $\M_{II}$ and $\M_{III}$. We need the following
\begin{oss}\label{segreconic}
We view the natural map
$$
(\det V_1)^2  \rightarrow S^2(S^2(V_1)) {\rm  \ given  \ by  \ }
(x_0 \wedge x_1)^2  \rightarrow   x_0^2 x_1^2 - (x_0x_1)^2
$$
(because $V_1 \cong E_{[0]}(2,1)$) as a map
$$\hol_B(2[0]) \rightarrow (\hol_B^3 \oplus L_1 \oplus L_2 \oplus L_3) (2[0]),
$$
given through  a vector with $6$ entries,
whose last $3$ are necessarily zero.
Moreover each of the first  $3$ is nonzero, since the Segre image of
$\PP^1$ in
$\PP^2$ is a smooth conic, i.e.,  a conic of rank $3$.

In particular (up to automorphisms) we can assume this vector to be
the transpose of $(1,1,1,0,0,0)$.
\end{oss}
\begin{lem}\label{unozero}
If $V_2(-[0])$ splits as $W \oplus L$, with $L \otimes L \cong \hol_B$ and
$W$ a vector bundle of rank $2$, then
$$\tilde{\A}_6 \cong \left( S^3(W) \oplus \left( S^2(W) \otimes
L\right) \right)
([0]-2\tau).$$
\end{lem}
\Proof
Applying  definition \ref{F} in this case, we get that
$\tilde{\A}_6$ is the cokernel of an injective map
$$
(W \oplus L) ([0]-2\tau)
\hookrightarrow
\left( [S^3(W) \oplus (S^2(W) \otimes L)] \oplus \{W \oplus L \} \right)
      ([0]-2\tau).$$

Applying remark \ref{segreconic} to the definition of the exact sequence
\ref{exseqAeven} for $n=2$, one easily sees that the second component
of this map is the identity, thus the above sequence
splits.
\QED

We can now discuss the families $\M_{II}$ and $\M_{III}$ separately.
\begin{prop}\label{M2}
$\M_{II}$ is either empty or it has dimension $4$.
\end{prop}
\Proof
By lemma \ref{classV2}, in case $II$, $\xi$ varies in a $1$-parameter
family.

By lemma \ref{unozero} and by the cited result of Atiyah (\cite{atiyah})
$$\tilde{\A}_6 \cong E_{\tau_i}(2,1)([0]-\tau) \oplus
E_{\tau_i}(2,1)([0]-\tau) \oplus
\left( (\oplus_{j=1}^3L_j) ([0]-\tau) \right),$$
and we obtain $h^0(\tilde{\A}_6)=2$ unless $\hol_{B}([0]-\tau)$ is a nontrivial
2-torsion bundle: in this last case
$h^0(\tilde{\A}_6)=3$.  Each of the two cases gives a (possibly
empty) unirational
family of dimension
$4$.
\QED

\begin{prop}\label{M3}
$\M_{III}$ has two connected components. Each is non-empty,
unirational of dimension $5$. For the first component  $\tau=[0]$,  whereas
for the second $\tau$ is a $2-$torsion point. In both cases
the branch curve $C = C_1 \cup C_2 \subset \CCC$ is disconnected.
\end{prop}
\Proof
By lemma \ref{classV2} we can assume $f_j=f_k=0$; applying lemma
\ref{unozero} we get
\begin{multline*}
\tilde{\A}_6 \cong
L_i([0]+\tau)
\oplus
L_j([0])
\oplus
L_i([0]-\tau)
\oplus
L_j([0]-2\tau)
\oplus \\
\oplus
L_k([0])
\oplus
\hol_B([0]-\tau)
\oplus
L_k([0]-2\tau)
.
\end{multline*}

Therefore either $h^0(\tilde{\A}_6)=4$ or $h^0(\tilde{\A}_6)=5$, the last case
occurring when $\tau=[0]$ or when $\hol_B([0]-\tau) \cong L_i$.

The decomposition $V_2(-[0])=L_i(\tau) \oplus L_j \oplus L_k$ yields natural
coordinates $y_i,y_j,y_k$ on $\PP(V_2)$. In these coordinates our
conic bundle has equation $f_0^2 y_i^2 + y_j^2 + y_k^2$: we note that
it has only one singular point (of type $A_1$).

Each line bundle summand of $S^3(V_2)$ corresponds to a monomial $y_iy_jy_k$.
Therefore the relative cubic given by the corresponding section of
$H^0(\tilde{\A}_6)$
is cut by a relative cubic on $\PP(V_2)$. If $h^0(\tilde{\A}_6)=4$ we
see that the only monomials allowed to have nonzero coefficient are
the monomials $y_i^3, y_i^2y_j,  y_i^2y_k$: in particular the
branch curve is not reduced (it contains $\{y_i=0\}$ twice) and the
generality conditions in definition \ref{adm5ple} are never fulfilled:
this case does not occur.

Therefore $h^0(\tilde{\A}_6)=5$ and either $\tau=[0]$ or
$\hol_B(\tau-[0])=L_i$.

If $\hol_B(\tau-[0])=L_i$ we have a relative cubic of the form
$ay_i^3+by_i^2y_j+cy_i^2y_k+dy_iy_j^2$, where $a$, $b$, $c$ and $d$
are global sections of the respective line bundles.

This cubic is clearly reducible and we write the corresponding curve
on the conic bundle as $C_1 \cup C_2$ where $C_1=\{y_i=0\}$ and
$C_2=\{ay_i^2+by_iy_j+cy_iy_k
+dy_j^2=0\}$.  
 Here $d$ is a
section of the trivial bundle and it is nonzero (or we would get the same contradiction as
in the previous case): whence
$C_1 \cap C_2 = \emptyset$.

$C_1$ is obviously smooth, while $C_2$ varies in a linear system that
has no fixed points on the conic bundle, therefore by Bertini's
theorem its general element is smooth. If we moreover assume that $a$
(section of $\hol_B(2[0])$) does not vanish in $[0]$, we ensure that
the curve does not pass through the singular point of the conic
bundle: therefore there exists an element in $H^0(\tilde{\A}_6)$
fulfilling the open
conditions listed in definition \ref{adm5ple}.

If $[0]=\tau$ we have the relative cubic of equation
$ay_i^3+by_i^2y_j+cy_i^2y_k+dy_iy_jy_k$: the proof of the existence of
a 'good' cubic is identical to  the previous one and again the branch curve is
disconnected.

The computation of the "number of parameters" is easy: $1$ parameter
for $B$, no
parameters for $\tau$ and $\xi$, $4$ from $H^0(\tilde{\A}_6)$.
\QED

\begin{oss}\label{reducible}
The example constructed in \cite{cc1} belongs to the second 
connected component of $M_{III}$($\tau \neq [0]$). The fact
that  the two singular points of
the branch curve in $B^{(2)}$, blown up by the rational map $\PP(V_1) \dashrightarrow
\CCC$, are contained in one of the curves $T_{\eta}$ is equivalent to the
fact that the extension class
$\xi$  is as in case III of lemma \ref{classV2}.

This holds true because the surface is  a minimal resolution
of the singularities of a double cover of $B^{(2)}$ branched on
$$
T_{\eta}+D_{t}+D_{t+\eta}+D_{t+\eta'}+D_{t+\eta+\eta'}+E_{2t+\eta},
$$
for some $t \in B$ and a pair of distinct non trivial $2-$torsion
elements $\eta, \eta'$. This curve has in fact two quadruple points at
the complete intersection $T_{\eta}\cap E_{2t+\eta}=
\{t,t+\eta\} \cup \{t+\eta',t+\eta'+\eta\}$.

The above mentioned branch curve, after removing the fibre
$E_{2t+\eta}$, stays in $|6D_0-E_{p}|$ for $p=\eta-4t$, and it is
therefore given by a map $\hol_B(p) \rightarrow S^6(E_{[0]}(2,1))$;
so the branch curve in $\CCC$ comes from a map
$\hol_B(p+3\tau)\cong L_i(3[0]+\tau) \rightarrow \A_6$.
Since $\A_6$ splits (following lemma \ref{unozero}) as
$\left( S^3(L_i (\tau)\oplus L_j) \oplus \left(
S^2(L_i (\tau) \oplus L_j) \otimes L_k \right) \right)
(3[0])$, we find $L_i(3[0]+\tau)$ as a direct summand in the left summand,
and this characterises the second case of proposition \ref{M3}.
\end{oss}
We note now that theorem
\ref{mainstream} shows that $\M_I^{\circ} \neq
\emptyset$: in fact it shows that there is a component of the moduli
space where the degree $0$ line bundle $\det V_1 (-\tau)$
can be chosen general, whereas
for the stratum $\M_{I,3}$ we have that
$\det V_1(-\tau)$ is a $3-$torsion line bundle, and we have just
shown that for each nonempty irreducible component contained in
$\M_{II} \cup \M_{III}$, $\det V_1(-\tau)$ is a $2-$torsion line
bundle. In the next proposition we analyse this stratum.
\begin{prop}\label{generalcase}
$\M_{I}^{\circ}$ forms an irreducible unirational family of dimension $5$.
\end{prop}
\Proof By lemma \ref{classV2} all $f_i$'s are different from zero, and
   $V_2(-[0]) \cong E_{\tau}(3,1)$.

By theorem 14' of \cite{atiyah} we have that
$$ S^2(E_{\tau}(3,1)) \oplus \Lambda^2 E_{\tau}(3,1)= E_{\tau}(3,1)
\otimes E_{\tau}(3,1) \cong E_{\tau}(3,2)^{\oplus 3},
$$
whence $S^2(E_{\tau}(3,1)) \cong E_{\tau}(3,2)^{\oplus 2}$,
$\Lambda^2(E_{\tau}(3,1)) \cong E_{\tau}(3,2)$, and that
$$ E_{\tau}(3,1) \otimes E_{\tau}(3,2) \cong \hol_B(\tau) \oplus
(\oplus_1^8 M_i(\tau)),
$$ where the $M_i$'s are the line bundles with $M_i^3 \cong \hol_B$, $M_i
\not\cong \hol_B$.

Consider the Eagon Northcott exact sequence
\begin{multline*}
0 \rightarrow \det E_{\tau}(3,1) \rightarrow \Lambda^2 E_{\tau}(3,1)
\otimes E_{\tau}(3,1) \rightarrow\\ \rightarrow E_{\tau}(3,1) \otimes
S^2(E_{\tau}(3,1)) \rightarrow S^3(E_{\tau}(3,1)) \rightarrow 0.
\end{multline*}

We have shown that the first $3$ bundles are direct sums of line
bundles of degree $1$: by cancellation we get
$$ S^3(E_{\tau}(3,1))= ( \oplus_1^2  \hol_B (\tau) )
\oplus \left( \oplus_1^8 M_i(\tau)
\right).
$$
By definition \ref{F} we have an exact sequence
\begin{multline}\label{exseqMi}
   0 \rightarrow E_{\tau}(3,1) \otimes \hol_B([0]-2\tau) \rightarrow \\
\rightarrow
\left( ( \oplus_1^2  \hol_B  ) \oplus (\oplus_1^8 M_i) \right) ([0]-\tau)
\rightarrow \tilde{\A}_6
\rightarrow 0.
\end{multline}

Since $3\tau \neq [0]$, $H^0(\tilde{\A}_6) \cong H^1(E_{\tau}(3,1) \otimes
\hol_B([0]-2\tau)) \cong \CC^2$.

Summing up, we can parametrize
$\M_I^{\circ}$ via a rational $5$-dimensional family
($1$ parameter for $w$, $2$ for
$\xi$, $1$ for
$\tau$ and $1$ for $B$).
\QED

We are left with the stratum  $\M_{I,3}$, computationally
more complicated. We will show that it has dimension  $ \leq 4$.

We will need the following algebraic lemma.
\begin{lem}\label{matriciona}
Let $B$ be a smooth elliptic curve, $\tau\in B$ a point,
$V_1=E_0(2,1)$. Fix $f_0 \in H^0(\hol_B(\tau))
\setminus \{0\}$, and $\forall i \in \{1,2,3\}$ $f_i\in H^0(\oplus_1^3
L_i(\tau))$.

Let then
$$
0
\rightarrow
S^2(V_1)
\stackrel{\sigma_2}{\rightarrow}
V_2
\rightarrow
\hol_{\tau}
\rightarrow
0.
$$
be the extension associated to $(f_0,f_1,f_2,f_3)$,
and let $\tilde{\A}_6$ be the vector bundle associated to $\sigma_2$ (as in
definition \ref{F}).

Then $\tilde{\A}_6$ is isomorphic to the cokernel of the map $F$ from
$$
\left( \hol_B \oplus L_1 \oplus
L_2 \oplus L_3 \oplus \hol_B \oplus L_3 \oplus L_2 \oplus \hol_B
\oplus L_1 \right) ([0]-3\tau)
$$
to
\begin{multline*}
( \hol_B \oplus L_1 \oplus L_2 \oplus L_3 \oplus \hol_B \oplus L_3
\oplus L_2 \oplus \hol_B \oplus L_1 \oplus \\ \oplus L_2 \oplus L_3
\oplus L_1 \oplus \hol_B \oplus L_1 \oplus L_3 \oplus L_2 )
([0]-2\tau)
\end{multline*}
given by the transpose of the matrix
$$M=\left(
\begin{array}{cccccccccccccccc}
\scriptstyle f_0&\scriptstyle f_1&\scriptstyle f_2&\scriptstyle
f_3&&&&&&&&&&&&\\ &\scriptstyle f_0&&&\scriptstyle f_1&\scriptstyle
f_2&\scriptstyle f_3&&&&&&&&&\\ &&\scriptstyle f_0&&&\scriptstyle
f_1&&\scriptstyle f_2&\scriptstyle f_3&&&&&&&\\ &&&\scriptstyle
f_0&\scriptstyle -f_3&&\scriptstyle f_1&\scriptstyle -f_3&\scriptstyle
f_2&&&&&&&\\ &&&&\scriptstyle f_0&&&&&\scriptstyle f_2&\scriptstyle
f_3&\scriptstyle -f_1&&\scriptstyle -f_1&&\\ &&&&&\scriptstyle
f_0&&&&\scriptstyle f_1&&\scriptstyle f_2&\scriptstyle f_3&&&\\
&&&&&&\scriptstyle f_0&&&&\scriptstyle f_1&&\scriptstyle
f_2&\scriptstyle f_3&&\\ &&&&&&&\scriptstyle f_0&&\scriptstyle
-f_2&&\scriptstyle f_1&&&\scriptstyle f_3&\scriptstyle -f_2\\
&&&&&&&&\scriptstyle f_0&&&&\scriptstyle f_1&&\scriptstyle
f_2&\scriptstyle f_3
\end{array}
\right).
$$

In particular $H^0(\tilde{\A}_6)$ is isomorphic to the kernel of $H^1 (F)$.
\end{lem}
\Proof To keep all the diagrams on the page, we
denote by ${\B}$ the vector bundle $End(V_1)=\hol_B \oplus L_1 \oplus
L_2 \oplus L_3$.

Our by now familiar exact sequence
$$ 0 \rightarrow \hol_B([0]-\tau) \rightarrow {\B}([0])
\rightarrow V_2 \rightarrow 0
$$ induces the exact sequence
$$ 0 \rightarrow S^2({\B})(3[0]-\tau) \rightarrow
S^3({\B})(3[0]) \rightarrow S^3(V_2) \rightarrow 0.
$$

We have the following diagram with exact rows and columns
(cf. definition \ref{F}) 
\begin{equation}\label{diagramma}
\begin{array}{ccccccccc}
&&&&&&0&&\\ &&&&&&\downarrow&&\\ 0&\rightarrow&
\hol_B(3[0]-\tau)&\rightarrow&{\B}(3[0])&\rightarrow&\det(V_1)^2
\otimes V_2&\rightarrow&0\\
&&&&&&\downarrow i_3&&\\
0&\rightarrow&S^2({\B})(3[0]-\tau)&\rightarrow&S^3({\B})(3[0])&\rightarrow&S^3(V_2)&\rightarrow&0\\ 
&&&&&&\downarrow&&\\
&&&&&&{\tilde{\A}_6}(2[0]+2\tau)&&\\ &&&&&&\downarrow&&\\ &&&&&&0&&
\end{array}.
\end{equation}

Note that all the vector bundles outside the last column are
direct sum of line bundles. We now give  explicit 'vertical' maps
$$ \alpha:
\hol_B(3[0]-\tau)
\rightarrow
S^2({\B})(3[0]-\tau),\ \ \
\beta:
{\B}(3[0])
\rightarrow
S^3({\B})(3[0])
$$ to be inserted in   (\ref{diagramma}) to enlarge it
to a bigger commutative diagram.

Recall that $i_n$ is defined by the formula $i_n((x_0 \wedge
x_1)^{\otimes 2} \otimes
q)=(\sigma_2(x_0^2)\sigma_2(x_1^2)-\sigma_2(x_0x_1)^2)q$, for $q \in
S^{n-2}(V_2)$. In particular  $i_2$  factors as
$$
(\det V_1)^2 \rightarrow S^2(S^2(V_1)) \rightarrow S^2(V_2).
$$

Remark \ref{segreconic} gave us an explicit form of the above 
(injective) map on the left
between
$\det V_1^2 \cong \hol_B(2[0])$, and
$$
S^2(S^2(V_1)) \cong S^2((L_1
\oplus L_2 \oplus L_3)([0])) \cong
(\hol \oplus L_3 \oplus L_2 \oplus \hol \oplus L_1 \oplus \hol)(2[0]).
$$
It has the matrix
$^t(1,0,0,1,0,1)$ (it is different from the previous  one  
because now  we use 
  the lexicographic
order).

We now define $\beta$ as the composition
of the natural maps
\begin{multline*}
(\det V_1)^2 \otimes \B([0])
\rightarrow
S^2(S^2(V_1)) \otimes \B([0])
\rightarrow\\
\rightarrow
S^2(\B([0])) \otimes \B([0])
\rightarrow
S^3(\B([0]))
\end{multline*}
where the first map is the one we have just described (tensored by $\B([0])$),
and the second map
is induced by the inclusion $S^2(V_1) \hookrightarrow
\hol_B([0]) \oplus S^2(V_1)=\B([0])$.

A similar splitting can be done for $\alpha$:
using again the lexicographic order we get
$$
S^2(\B)(3[0]-\tau)=
(\hol_B \oplus L_1 \oplus L_2 \oplus L_3 \oplus \hol_B \oplus L_3
\oplus L_2 \oplus \hol_B \oplus L_1 \oplus \hol_B)(3[0]-\tau)
$$
and $\alpha = \ ^t ( 0, 0, 0, 0, 1, 0, 0, 1, 0, 1).$

Taking the mapping cone we obtain a free resolution of $\tilde{\A}_6$ 
as follows:

\begin{multline*}
0
\rightarrow
\hol_B(3[0]-\tau)
\stackrel{\tilde{\alpha}}{\rightarrow}
\B([3[0]) \oplus S^2(\B)(3[0]-\tau)
\rightarrow\\
\stackrel{\tilde{\beta}}{\rightarrow}
S^3(\B)(3[0])
\rightarrow
\tilde{\A}_6(2[0]+2\tau)
\rightarrow
0
\end{multline*}
with
$$ \tilde{\alpha}= \begin{pmatrix} f_0\\ f_1\\ f_2\\ f_3\\
0\\ 0\\ 0\\ 0\\ -1\\ 0\\ 0\\ -1\\ 0\\ -1
\end{pmatrix}
\tilde{\beta}= \left(
\begin{array}{cccccccccccccc}
   & & & &f_0&   &   &   &   &   &    &   &   &   \\
   & & & &f_1&f_0&   &   &   &   &    &   &   &   \\
   & & & &f_2&   &f_0&   &   &   &    &   &   &   \\
   & & & &f_3&   &   &f_0&   &   &    &   &   &   \\
1& & & &   &f_1&   &   &f_0&   &    &   &   &   \\
   & & & &   &f_2&f_1&   &   &f_0&    &   &   &   \\
   & & & &   &f_3&   &f_1&   &   &f_0 &   &   &   \\
1& & & &   &   &f_2&   &   &   &    &f_0&   &   \\
   & & & &   &   &f_3&f_2&   &   &    &   &f_0&   \\
1& & & &   &   &   &f_3&   &   &    &   &   &f_0\\
   &1& & &   &   &   &   &f_1&   &    &   &   &   \\
   & &1& &   &   &   &   &f_2&f_1&    &   &   &   \\
   & & &1&   &   &   &   &f_3&   &f_1 &   &   &   \\
   &1& & &   &   &   &   &   &f_2&    &f_1&   &   \\
   & & & &   &   &   &   &   &f_3&f_2 &   &f_1&   \\
   &1& & &   &   &   &   &   &   &f_3 &   &   &f_1\\
   & &1& &   &   &   &   &   &   &    &f_2&   &   \\
   & & &1&   &   &   &   &   &   &    &f_3&f_2&   \\
   & &1& &   &   &   &   &   &   &    &   &f_3&f_2\\
   & & &1&   &   &   &   &   &   &    &   &   &f_3
\end{array}
\right) .
$$

The reader can now easily simplify this resolution and get the minimal one that
is given by the above matrix $^t M$.
\QED

We now consider the missing case  $\M_{I,3}$:
$V_2$ is indecomposable, and $\tau$ is a 3-torsion
point. Exact sequence (\ref{exseqMi}) still holds
but the vector bundle in the middle has nontrivial cohomology:
both cohomology groups have the same dimension: $1$  if $\tau \neq
[0]$, and $2$ if $\tau =[0]$.
\begin{prop}\label{tau=0}
  $\M_{I,3}$ has dimension at most $4$.
\end{prop}
\Proof
We start with the case where $\tau$ is a 3-torsion point, $\tau 
\neq [0]$ and assume by contradiction
that we have then an irreducible family $\SSS$ of dimension at least $5$.
The data $(B,\tau,\xi)$ give a map from
$\SSS$ to an irreducible family $\mathcal Y$ of dimension $3$ (1-parameter for
$B$, no parameters for $\tau$, $2$ parameters for $\xi$) whose fibre
is the projective space  $\PP(H^0(\tilde{\A}_6))$ which has dimension $2$.

We know in fact that $h^0(\tilde{\A}_6) \leq 3$, that the map 
dominates $\mathcal Y$,
whence the fibre has dimension exactly two.

In other terms, for any triple $(B,\tau,\xi)$ with $\tau\neq [0]$ of 
$3-$torsion we have $h^0(\tilde{\A}_6) = 3$.

Fix then a general pair ($B$, $\tau$) with $\tau$ a
nontrivial $3-$torsion
point. In lemma \ref{matriciona} we have shown that $h^0(\tilde{\A}_6)$ equals
the dimension of the kernel of $H^1(F)$, where  $F$ is determined by
$f_0,f_1,f_2,f_3$: therefore we are assuming that for our choices of 
$B,\tau$ and general
choice of the $f_i$'s $H^1(F)$ has a kernel of dimension
$3$. Letting now $f_1$ tend to zero we obtain, in the limit,   case
$\M_{II}$  where, as shown in the proof of proposition
\ref{M2}  (since in this case $\tau$ is not a $2-$torsion),
$h^0(\tilde{\A}_6)=2$: this contradicts the semicontinuity theorem.

The case $\tau=[0]$ is more difficult, since in this case we
only know  that $2 \leq h^0(\tilde{\A}_6)\leq 4$.

Also in this case we use lemma \ref{matriciona}.
Dualizing, we see that $H^0(\tilde{\A}_6)$ is
isomorphic to the cokernel of the map $F'$ from
\begin{multline*}
H^0(( \hol_B \oplus L_1 \oplus L_2 \oplus L_3 \oplus \hol_B \oplus L_3
\oplus L_2 \oplus \hol_B \oplus L_1 \oplus \\ \oplus L_2 \oplus L_3
\oplus L_1 \oplus \hol_B \oplus L_1 \oplus L_3 \oplus L_2 )
([0]))
\end{multline*}
to
$$
H^0\left( \hol_B \oplus L_1 \oplus
L_2 \oplus L_3 \oplus \hol_B \oplus L_3 \oplus L_2 \oplus \hol_B
\oplus L_1 \right) (2[0])
$$
given by the matrix $M$.

Since we are assuming $V_2$ to be indecomposable, all the $f_i$'s are not zero
and therefore  each $f_i$ generates the corresponding $H^0$.
We can then write explicitly generators of the image of $F'$ as vectors whose
entries have the form $f_if_j$.

The pairs $\{f_0f_i,f_jf_k\}$ give bases of $H^0(L_i(2[0]))$ (for
$\{i,j,k\}=\{1,2,3\}$); the remaining $4$ elements
($f_0^2,f_1^2,f_2^2,f_3^2$) generate $H^0(\hol_B(2[0]))$: therefore
there are two independent relations among them.

We identify $(f_0,f_1,f_2,f_3)$ to scalar
multiples of the functions
$$(\theta_{11},\theta_{10},\theta_{01},\theta_{00}),$$ where
the $\theta_{ij}$'s  are the half-integer theta functions
$\theta_{ij}(z,\mu)$ (cf. e.g. \cite{tatamum}, page 17: here $\mu$
is a point of  Poincar\'e 's  upper half plane).

Therefore the relations between the $f_i$'s come from the relations
between the $\theta_{ij}$'s; using the relations $(E_1)$ and $(E_2)$ 
in \cite{tatamum}
page 23, we can write them as  $f_2^2=af_0^2+bf_1^2$,
$f_3^2=cf_0^2+df_1^2$, with $a,b,c,d$ that vary freely in an open set
of $\C^4$ according to the choice of $B$ and  of the $f_i$'s.

We have now given a basis of each respective space of global sections, and
we may write explicitly (depending on the parameters $a,b,c,d$) the
corresponding matrix of $F'$.
We wrote a Macaulay 2 script (available upon request) that told us
\begin{itemize}
\item
that $F'$ is generically injective ($h^O({\tilde{\A}_6})=2$),
\item
$F'$ acquires a kernel of dimension $1$ ($h^O({\tilde{\A}_6})=3$)
on a subvariety   of codimension $1$ of the space of parameters
\item
$F'$ acquires a kernel of dimension $2$ in codimension $4$
(therefore never, the parameters $a,b,c,d$ giving only three moduli up to
automorphisms).
\end{itemize}

We get then two (possibly empty) families:

$i) h^O(\tilde{\A}_6)=2$ gives a family of dimension at most $3+1=4$;

$ii) h^O(\tilde{\A}_6)=3$ gives a family of dimension at most $3-1+2=4$.
\QED

We can now conclude our classification theorem

{\it End of the proof of theorem \ref{Kquadro3gen2}.}
We have shown the existence of three irreducible components of dimension $5$:
the main stream component whose general point is contained in $\M_I^{\circ}$,
and the two ones contained in the closed set $\M_{III}$, according to proposition
\ref{M3}.

We have also seen that every other family in the moduli space has 
dimension at most
$4$; by remark \ref {strategy} we conclude that $\M'$ has three
irreducible components of dimension $5$.

It remains to show that they do not intersect.
The two components in $\M_{III}$  do not intersect,
since (proposition \ref{M3}) in one case $\tau$ is $[0]$, in the other
case it is a non trivial  $2-$torsion point.

Since a general point of the main stream component is contained in
$\M_I^{\circ}$ it is then enough to show that
$\overline{\M_I^{\circ}} \cap \M_{III}=\emptyset$. It will be then enough
to show that there is no flat family on a disc whose central fibre is a
surface in $\M_{III}$ and whose general fibre is a general surface in
$M_I^{\circ}$.

We borrow an argument used often in  Horikawa's work.
We note that by theorem
\ref{mainstream} a general surface of type $M_I^{\circ}$ is a double
cover of $\CCC$ with irreducible branch curve. On the contrary we
have shown in  proposition \ref{M3} 
that the same double cover, for a surface in
$\M_{III}$, has a disconnected branch curve $C_1 \cup C_2$. Therefore we
would have a family of double covers with a connected general branch curve
  and with a disconnected special branch curve, a contradiction.
\qed

\section{Fibrations with non hyperelliptic general fibre:
      the case $g=3$}
In the rest of the paper we will concentrate on the case  when
  the general fibre of $f:S \rightarrow B$ is nonhyperelliptic of
  genus $g=3$.

The natural morphism of graded $\hol_B$-algebras $\sigma:
    Sym(V_1) \rightarrow {\R}(f)$ has as kernel the graded sheaf of ideals
  $\L$ and we denote as usual by
        $\T$ its (graded) cokernel.
\begin{oss}\label{notazgen3}
Letting $ X : ={ \bf Proj}\ (\R(f))$  the rational map
\[\psi_1: X \dashrightarrow \Sigma \subset
{\bf Proj}(Sym\ V_1)\ =\PP(V_1)\]
allows to factor
the relative canonical map $\varphi : S \dashrightarrow \Sigma$ as
$\psi_1 \circ r$, where  $r: S \ra X$ is a minimal
resolution of singularities (RDP's ).

Observe that, if the general fibre of $f$ is nonhyperelliptic, all the
above maps are birational and  $\L$ is the ideal sheaf of the 
canonical image $\Sigma:=
\varphi(S)$.
\end{oss}

The next lemma  investigates the sheaf $\T$: assuming the general
fibre to be nonhyperelliptic, we know that the $\T_n$ are all torsion
sheaves on $B$. 

\

{\bf Main assumption:} We will assume (often without explicit mention) in
the rest of the section that every fibre of the genus $3$ fibration
$f$ is $2-$ connected.

\

\begin{lem}
\label{contact2}
Let $f$ be a genus 3 fibration whose general fibre is nonhyperelliptic and
such that every fibre is 2-connected. Then:

1) the sheaf $\T_2$ is isomorphic to the structure sheaf $\hol_{\tau}$ of an
   effective  divisor $\tau$ on $B$;

2) the sheaves $\T_n$ are isomorphic to free
$\hol_{\tau}$-modules of rank $2n-3$;

3) the multiplication map $V_1 \otimes V_2 \rightarrow V_3$ induces an
       isomorphism $V_1 \otimes \T_2 \cong \T_3$.
\end{lem}
\Proof The argument is similar to the one given in Lemma \ref{torsgenus2}.

By the classification of genus $3$ fibres due to M. Mendes Lopes
(cf. \cite{mml}), and by the hypothesis of
2-connectedness, a fibre $F$ is either nonhyperelliptic, i.e., it
has a canonical ring of the form
\[
\C[x_1,x_2,x_3]/<F_4(x_i)>,
\]
or it is honestly hyperelliptic, i.e.  it has a canonical ring of the form
\[
R=\C[x_1,x_2,x_3,y]/<r_1:=Q(x_i),r_2:=y^2-G(x_i)>,
\]
where $\deg x_i=1$, $\deg y=2$, $\deg Q=2$, $\deg G=4$.

If $t$ is a local parameter in $B$ such that the point $p:=\{t=0\}$ is
the image of a hyperelliptic fibre, the relation $r_1$ lifts to a
relation
\[
\overline{r}_1=\overline{Q}(x_i,t)+\mu(t)y
\]
where $\overline{Q}(x_i,0)=Q(x_i)$, and $\mu(0)=0$. The
assumption that the generic fibre is nonhyperelliptic imposes $\mu
\not\equiv 0$, whence we may assume as in \ref{torsgenus2} $ \mu(t) = t^s$
and then $(\T_2)_p=\hol_{B,P}/t^s$, proving the first part of
the statement.

Using the lift of $r_2$ to eliminate the multiples of
$y^2$,  we see that, chosen  a basis $\{q_{j}\}$
for the homogeneous part of degree $k-2$
of the quotient ring $\C[x_1,x_2,x_3]/(Q)$,
the set $\{t^iq_{j}y | i < s \}$ is a
basis for the complex vector space, stalk of $\T_k$ at $p$. It follows
  then that ${\T_k}_{,p}$ is a free $\hol_{B,p}$-module with basis
$\{q_{j}y\}$.

For $k=3$ the above basis is $\{x_1y,x_2y,x_3y\}$, i.e.,  exactly
the image of the natural basis of $V_1 \otimes \T_2$, and
the lemma is proven.
\qed

The next lemma  investigates the sheaf $\L$.  Write
$\pi:\PP(V_1) \rightarrow B$ for the canonical projection.

\begin{lem}\label{kergen3}
If $f$ is a genus $3$ fibration with
nonhyperelliptic general fibre and such that every fibre is
2-connected, then

1) $\L_2=\L_3=0$.

2) $\deg V_1=\chi-2(b-1)$, $\deg \tau=K^2-3\chi(\hol_S)-10(b-1)$.

3) $\L_4$ is a line bundle of degree $\deg V_1 - \deg \tau$.

4) the maps $\L_{4} \otimes S^{n-4}(V_1) \hookrightarrow
\L_n$ are isomorphisms: in particular the ``relative
canonical image'' $\Sigma$ is a divisor on ${\PP} (V_1)$
belonging to the linear system $|\hol_{\PP(V_1)}(4) \otimes \pi^*
\L_4^{\textsc v}|$.
\end{lem}
\Proof 
We have already noted that all the maps $\sigma_n$ are
generically surjective, so that their kernels are vector bundles of
rank  $  { n+2 \choose 2} - 4n + 2 $ for $n \geq 2$. We
get
\begin{itemize}
\item $\L_2=\L_3=0$;
\item $\L_4$ is a line bundle;
\item $\forall n \geq 4$, rank $\L_n=$rank $S^{n-4}(V_1)$.
\end{itemize}

By remark \ref{fujitadeg} $\deg V_1 = \chi(\hol_S) -2(b-1),\ \deg
V_2=\chi(\hol_S)+K^2-18(b-1)$ and once more
the exact sequence
$ 0 \rightarrow S^2(V_1) \rightarrow V_2 \rightarrow \hol_{\tau}
\rightarrow 0
$ yields
$ \deg \tau =\deg V_2-4\deg V_1 = K^2-3\chi(\hol_S)-10(b-1).$

By remark \ref{fujitadeg} $\deg
V_2=6K^2+\chi(\hol_S)-98(b-1)$. Moreover $\deg S^4(V_1)=20 \deg V_1 =
20\chi(\hol_S)-40(b-1)$. By lemma \ref{contact2} $\deg \T_4=5
\deg{\tau}=5K^2-15 \chi(\hol_S)-50(b-1)$, and we conclude
\begin{multline*}
\deg \L_4=\deg S^4(V_1)+\deg \T_4-\deg V_4=\\
=5(\chi(\hol_S)+K^2)-90(b-1)-6K^2-\chi(\hol_S)+98(b-1)=\\
=4\chi(\hol_S)-K_S^2+8(b-1)=\deg V_1-\deg \tau.
\end{multline*}

4) is obvious since $\Sigma$ is a Cartier divisor in $\PP (V_1)$.
\qed

Now, with the informations on $\T$ and $\L$
provided by lemmas \ref{contact2} and lemma \ref{kergen3}, we can
investigate the $\hol_B$-algebra structure of $\R(f)$, i.e the
multiplication maps $\mu_{i,j}:V_i \otimes V_j \rightarrow V_{i+j}$.

For $i=j=1$ we have
the composition $V_1 \otimes V_1 \rightarrow S^2(V_1)
\stackrel{\sigma_2}{\rightarrow} V_2$.

The next proposition shows that the triple ($V_1,V_2,\sigma_2$) determines
$\R(f)$ in degree $\leq 3$.

\begin{df}
Let $A: V_1 \otimes \Lambda^2(V_1) \rightarrow
S^2(V_1) \otimes V_1$ be defined by
$$ A(c \otimes (a \wedge b)): = bc\otimes a - ac \otimes b
$$ and $B: \Lambda^3(V_1) \rightarrow V_1 \otimes \Lambda^2(V_1)$
be defined by
$$ B(a \wedge b \wedge c):= a \otimes (b \wedge c) + b \otimes (c
\wedge a) + c \otimes (a \wedge b).
$$
\end{df}

\begin{prop}\label{deg3gen3}
$V_3=$coker $((\sigma_2 \otimes Id) \circ A):V_1 \otimes
\Lambda^2(V_1) \rightarrow V_2 \otimes V_1)$ and $\mu_{2,1}$ is given by the
projection onto the cokernel.
\end{prop}
\Proof Consider the diagram
$$\begin{array}{ccccccccc} &&0&&0&&0&&\\
&&\uparrow&&\uparrow&&\uparrow&&\\ 0&\rightarrow
&S^3(V_1)&\stackrel{\sigma_3}{\longrightarrow}&V_3&\longrightarrow&{\T}_3&
\rightarrow &0\\ &&\uparrow&&\ \ \ \ \ \uparrow {\scriptstyle
\mu_{2,1}}&&\uparrow&&\\ 0&\rightarrow&S^2(V_1) \otimes
V_1&\stackrel{\sigma_2 \otimes Id}{\longrightarrow}&V_2 \otimes
V_1&\longrightarrow&\T_2 \otimes V_1&\rightarrow &0\\ &&\ \
\ \uparrow {\scriptstyle A}&&\ \ \ \ \ \ \ \ \ \ \ \ \uparrow
{\scriptstyle (\sigma_2 \otimes Id) \circ A}&&\uparrow&&\\
0&\rightarrow &V_1 \otimes \Lambda^2(V_1)& = &V_1 \otimes
\Lambda^2(V_1)&\longrightarrow&0& &\\ &&\ \ \ \uparrow{\scriptstyle
B}&&\uparrow&&&&\\ && \Lambda^3(V_1) &=&\Lambda^3(V_1)&&&&\\
&&\uparrow&&\uparrow&&&&\\ &&0&&0&&&&.
\end{array}
$$ where the map $S^2(V_1) \otimes V_1 \rightarrow S^3(V_1)$ is the
natural one in $Sym(V_1)$ while the map $\T_2 \otimes V_1
\rightarrow \T_3$ is the push forward of $\mu_{2,1}$
(it is an isomorphism by lemma \ref{contact2}).

It is then easy to check that the diagram commutes with exact rows and
columns. This proves our statement.
\qed

We  now analyse the algebra $\R(f)$ in degree $4$.
\begin{df}
Let $C:S^2(\Lambda^2 (V_1)) \rightarrow
S^2(S^2(V_1))$ be defined by
$$C((a \wedge b)(c \wedge d)) := (ac)(bd)-(ad)(bc).$$
\end{df}
\begin{lem}\label{Flocfree}
The map $S^2(\sigma_2) \circ C: S^2(\Lambda^2 (V_1)) \rightarrow
S^2(V_2)$ is injective with locally free cokernel.
\end{lem}
\Proof We must show that the map $C$ is injective on every fibre. If
$\phi$ is a map between vector bundles on $B$, for every $p \in B$ we
will denote with $(\phi)_{(p)}$ the corresponding linear map
between the fibres over $p$.

$C$ is  injective on every fibre because rank $V_1 = 3$. Note that
injectivity  fails for higher rank since, if $x_1,
\dots x_4$ are linearly independent, then
$$\sum_{\sigma (1)=1 , sign(\sigma) = 1}(x_{1} \wedge x_{\sigma(2)}) 
(x_{\sigma(3)}
\wedge x_{\sigma(4)})
\mapsto$$
$$ \mapsto \sum_{\sigma (1)=1 , sign(\sigma) = 1}(x_1  x_{\sigma(2)})
(x_{\sigma(3)} x_{\sigma(4)}) - (x_1  x_{\sigma(4)}) (x_{\sigma(2)} 
x_{\sigma(3)}) =
0.$$

We saw in the proof of lemma \ref{contact2} that, if $F_p$ is not
hyperelliptic, $(\sigma_2)_{(p)}$ is an isomorphism, so
$(S^2(\sigma_2) \circ C)_{(p)}$ is injective.

At points $p$ of $B$ for which the fibre $F_p$ is hyperelliptic, the
kernel of $(\sigma_2)_{(p)}$ has rank $1$.
We can choose a basis $x_0,x_1,x_2$ for the fibre of $V_1$ over $p$
so that this kernel is generated either by $x_0^2$, or by
$x_0^2+x_1^2$, or by $x_0^2+x_1^2+x_2^2$, according to the rank of the
conic $q$ canonical image of the corresponding hyperelliptic fibre.

This choice induces a basis of $S^2(\Lambda^2(V_1))$ given by three
vectors of the form $(x_i \wedge x_j)(x_i \wedge x_j)$ and three of
the form $(x_i \wedge x_j)(x_i \wedge x_k)$.

By definition:
$$C((x_i \wedge x_j)(x_i \wedge x_j))=(x_i^2)(x_j^2)-(x_ix_j)(x_ix_j)$$
$$C((x_i \wedge x_j)(x_i \wedge x_k))=(x_i^2) (x_jx_k)-(x_ix_j)(x_ix_k)$$

The $6$ image vectors  of the basis elements are linearly 
independent.  It remains to show that  the  subspace they span
intersects transversally the subspace of $S^2(S^2(V_1))$ consisting of
multiples of $q$. This is  straightforward: it suffices to send to
zero the subspace spanned by the elements $(x_i^2) (x_j x_k)$ (for all
$i,j,k$). \qed

\begin{df}\label{Fgen3}
1) Given a curve $B$, a rank $3$ vector bundle $V_1$ on $B$, an
effective divisor $\tau$, and an extension
$$ 0 \rightarrow S^2(V_1) \rightarrow V_2 \rightarrow \hol_\tau
\rightarrow 0
$$
where $V_2$ is still a vector bundle, set $$\tilde{V_4} : =
S^2(V_2)/S^2(\Lambda_2(V_1))= \ { \rm coker }\ (S^2(\sigma_2) \circ
C).$$  Lemma \ref{Flocfree}  shows that $\tilde{V_4}$ is a locally
free vector bundle of rank $15$. 

2) Note that coker $C = S^4(V_1)$, so the image of $S^2(\sigma_2) \circ
C$ is in the kernel of the map $S^2(V_2) \rightarrow V_4$, inducing a
map $\tilde{V_4} \rightarrow V_4$: we  denote by $\L'_4$ its kernel.
\end{df}
\begin{prop}\label{deg4gen3}
\begin{itemize}
\item[1)] The maps $S^2(V_2) \rightarrow V_4$ and $\tilde{V_4} \rightarrow V_4$
    are surjective.
\item[2)] $\L_4 \cong \det V_1 \otimes \hol_B(-\tau)$.
\item[3)] $\L'_4 \cong \det V_1 \otimes \hol_B(\tau)$.
\item[4)] Consider the two linear maps induced, on the fibres  over a point $p \in Supp
(\tau)$,  by the embedding
$\L'_4 \hookrightarrow \tilde{V_4}$
    and by the composition map $V_2 \otimes S^2(V_1) \ra
    S^2(V_2) \ra \tilde{V_4}$. Their images are  vector subspaces 
intersecting only
    in zero.
\end{itemize}
\end{prop}
\Proof The first part is an immediate consequence of the fact that
$\forall p$ the map $S^2(H^0(F_p,\omega^2_{F_p})) \rightarrow
H^0(F_p,\omega^4_{F_p})$ is surjective, as made clear by the explicit
description of the canonical rings given in the proof of lemma
\ref{contact2}.

\

Recall that the map $\L_4 \hookrightarrow S^4(V_1)$ defines
$\Sigma$ as a divisor in $|\hol_{{\PP}(V_1)}(4) \otimes
\pi^* ({\L_4})^{\textsc v}|$, as shown in lemma \ref{kergen3}.

The dualizing sheaf of ${\PP}(V_1)$ is $\omega_{{\PP}(V_1)}=
\hol_{{\PP}(V_1)}(-3) \otimes \pi^* (\det V_1
\otimes \omega_B ) $ (cf. \cite{har}, ex. III.8.4.(b)), therefore the
dualizing sheaf of $\Sigma$ is the restriction of
$\hol_{{\PP}(V_1)}(1)\otimes \pi^* (\L_4^{\textsc
v}\otimes \det V_1 \otimes \omega_B)$.

The morphism (cf. remark \ref{notazgen3}) $\psi_1:X \rightarrow
\Sigma$ is an isomorphism when restricted to the preimage of $B
\setminus Supp(\tau)$, and 
since the relative dualizing sheaf  $\omega_{S|B}$ is
on the one side 
$$ \omega_{S|B}= \varphi^*\hol_{\PP}(1),$$ on the other side 
it equals 
$$ \mathcal I  \omega_{\Sigma | B} \cong \omega_{S|B} \cong \mathcal I 
\varphi^*
\hol_{\PP}(1)\otimes f^*((\L_4)^{\textsc v} \otimes \det V_1),
$$
where $\mathcal I$ is the conductor ideal,
2) follows if we prove that  $\mathcal I$ is the principal ideal
associated to the divisor $\tau$.

We will now use the local equations for $X$ in a
neighborhood of a fibre $F_p$ with  $ p \in Supp (\tau)$.

Choose a small open neighborhood $U$ of $p$ such that  
$X$ is the subvariety of $ U \times {\PP}(1,1,1,2)$
defined by ideal generated by $t^sy- \bar{Q}(x,t)$ 
and $y^2-\bar{G}(x,t)$. The above equations
show that $\hol_X$ is locally generated
by  $\{ 1, y\}$ as an $\hol_{\Sigma} $-module
  and the rows of the following matrix
are relations for $\{ 1, y\}$:
$$
\begin{pmatrix}
-t^s \bar{G}&\bar{Q}\\
-\bar{Q}&t^s
\end{pmatrix}.
$$

Since the determinant of the above matrix is a generator
of the ideal of $\Sigma$ in  $U \times {\PP}^2 $, 
it follows that the above is the full matrix
of relations, and we conclude that the conductor
ideal is generated by  $\bar{Q}, t^s $ as an ideal in $\hol_{\Sigma} $
and by  $t^s $ as an ideal in $\hol_X$.

\

We have the  commutative diagram with exact rows and columns:
$$\begin{array}{ccccccccccc} &&&&&&0&&0&&\\
&&&&&&\uparrow&&\uparrow&&\\ 0&\rightarrow & {\L}_4&\rightarrow
&S^4(V_1)&\stackrel{\sigma_4}{\longrightarrow}&V_4&\rightarrow&{\T}_4&\rightarrow 
&0\\ &&&&||&&\uparrow&&\uparrow&&\\
&&0&\rightarrow&\frac{S^2(S^2(V_1))}{S^2(\Lambda^2V_1)}&
\stackrel{S^2(\sigma_2)_{/S^2(\Lambda^2V_1)}}{\longrightarrow}&{\tilde{V_4}}&\rightarrow& 
{\K}&\rightarrow& 0\\ &&&&&&\uparrow&&&&\\
&&&&&&\L'_4&&&&\\ &&&&&&\uparrow&&&&\\ &&&&&&0&&&&
\end{array}
$$

By lemma \ref{contact2}
$ {\T}_4= {\hol_{\tau}}^{\oplus 5},
$ and
$\K$ is isomorphic to coker $S^2(\sigma_2)$, whence
$$ \K \cong {\hol_{\tau}}^{\oplus 5} \oplus {\hol_{2\tau}}.
$$

Computing determinants in the three exact sequences of the diagram
above, one finds:
$$\begin{array}{l} \det V_4 \otimes \L'_4= \det \tilde{V_4}= \det S^4(V_1)
     \otimes \hol_B(7\tau) \\ \det V_4 \otimes
\L_4= \det S^4(V_1) \otimes \hol_B(5\tau)
\end{array}$$
therefore $\L_4=\L'_4(-2\tau)$ and also 3) is proven.

\

Finally we want to prove that the embedding $\L'_4
\hookrightarrow \tilde{V_4}$ satisfies 4). 
Consider $p \in supp (\tau)$. The fibre over $p$ of
$\R(f)$, i.e., the stalk modulo the maximal ideal of $p$,
is the canonical ring of the corresponding fibre  $F_p$. In
degree $\leq 4$ we can now see that it has three generators
$x_1,x_2,x_3$ in degree $1$ given by $V_1$, one generator $y$ in degree
$2$,  generating $\T_2$, one relation, $Q(x_i)$,
generating the kernel of the map induced by $\sigma_2$ at the level of
the fibres, and one further relation in degree $4$, given by a
representative in $S^2(V_2)$ of the image of a generator of $\L'_4$.
We know by the 2-connectedness hypothesis that this relation has
the form $y^2=G(x_i)$, whence the image of $V_2 \otimes S^2(V_1)$ is
given by the polynomials without the term $y^2$.
\qed

\begin{df}\label{datagen3} 
Define the  {\bf genus three 5-tuple} associated to a
genus $3$ fibration with nonhyperelliptic general fibre and 
 2-connected fibres as ($B$, $V_1$, $\tau$,
$\xi$, $w$),where
\begin{itemize}
\item
$\xi$ is the element $\xi \in Ext^1(\hol_{\tau},S^2(V_1)) / Aut (\hol_{\tau})$
yielding the exact sequence $0 \rightarrow
S^2(V_1) \rightarrow f_*(\omega_{S|B}^2) \rightarrow \hol_{\tau}
\rightarrow 0$,
\item
$w$ is the element  in $\PP ( H^0 ( \tilde{V_4} \otimes 
{\L'}_4 ^{-1} ))$ corresponding 
to the embedding $i : \L'_4 \hookrightarrow
      \tilde{V_4} = {S^2(V_2)/S^2(\Lambda^2(V_1))}$. 
\end{itemize}
\end{df}

We have seen that
\subitem i) $\xi$ yields a vector bundle ($V_2$);
\subitem ii) $\L'_4 \cong \hol_B(\tau) \otimes \det V_1$;
\subitem iii) $i$ has a locally free cokernel;
\subitem iv) the image of   the fibre map corresponding to $i$
at a point in $ supp(\tau)$ is not contained in the image
of $V_2 \otimes S^2(V_1)$.

The next theorem  shows that these data determine the fibration.
\begin{df}\label{X3}
1) Define a {\bf  genus three 5-tuple}
($B$, $V_1$, $\tau$, $\xi$, $w$) as follows:
\begin{itemize}
\item $B$ is a smooth curve,
\item $V_1$ is a rank $3$ vector bundle,
\item $\tau$ is an effective divisor on $B$,
\item $\xi$ is an element of $ Ext^1(\hol_{\tau},S^2(V_1)) / Aut (\hol_{\tau})$
      yielding a vector bundle $V_2$,
\item if $\tilde{V_4}$ is the vector bundle constructed via $\sigma_2$ as
in definition \ref{Fgen3}, then $w = \PP(i)$, where $i$ is an 
embedding $\hol_B(\tau)
\otimes \det V_1 \hookrightarrow \tilde{V_4}$
such that iii) and iv) above hold.
\end{itemize}

2) Define its {\bf associated relative canonical model} $X$ as follows:
let $W \subset \PP (V_2)$ be the image of the rational Veronese map
$ \PP (V_1) \dashrightarrow \PP (V_2)$ induced by $\sigma_2$,  and let
  $X$ be the relative quadric divisor on $W$ cut by the 
principal ideal
generated by the image of $i$.

3) Define a   genus three 5-tuple to be  {\bf  admissible} if
its associated relative canonival model $X$ has
only Rational Double Points as singularities.
\end{df}
\begin{oss}
One can explicitly, cf. the proof of \ref{strucgen3},
  define a sheaf of graded $\hol_B$- algebras $\R$,
generated by $ Sym (V_1 \oplus V_2)$ such that $X = \Proj(\R)$. The 
geometric procedure
is more suitable to determine the singularities of $X$.
\end{oss}
\begin{teo}\label{strucgen3}
  Let $f$ be a relatively minimal genus $3$ nonhyperelliptic
fibration such that every fibre is $2-$connected.
Then its
{\bf associated 5-tuple} $(B,V_1,\tau,\xi,w)$ is admissible.

Viceversa, every admissible genus three 5-tuple
($B$, $V_1$, $\tau$, $\xi$, $w$)
is the associated 5-tuple of a unique  genus
$3$ nonhyperelliptic fibration
$f:S \rightarrow B$ with the property that every fibre is $2$-connected
and with invariants $\chi (\hol_S) =\deg V_1 +
2(b-1)$, $K^2_S=3\deg V_1 + \deg \tau +16 (b-1)$.
As in the case of genus 2, the bijection thus obtained is  functorial.
\end{teo}
\Proof
The 5-uple of such a fibration is admissible by definition.

Conversely, as in the proof of  theorem \ref{gen2}, we construct 
$V_2$ via the extension class $\xi$, 
$\tilde{V_4}$ via $\sigma_2$ as
in definition \ref{Fgen3}, 
and finally $X$ as
in definition \ref{X3}. 

By our assumption, $X$ has at most rational double
points as singularities and a minimal resolution $S$ of the
singularities of $X$ yields a genus $3$ fibration $f:S
\rightarrow B$ such that the general fibre is nonhyperelliptic . 

Let us show that each fibre $F$ of $f$ is  $2-$connected: 
otherwise one of following holds:
\begin{itemize}
\item
0) $F$ is not $1$-connected, i.e., by Zariski's lemma
$ F = 2 F'$, where $F'$ has arithmetic genus $2$ or
\item
1) $ F = A + B $, $ A \cdot B = 1$.
\end{itemize}

In case 0), since $ \hol_{F'} (2 K_S ) =  \hol_{F'} (2 K_{F'} )$,
the relative bicanonical map is not birational, a contradiction.

In case 1), the relative canonical map has base points,  contradicting 
the fact that, by iv), the coefficient of $y^2$ in the local equations
is equal to $1$.

The rest of the proof is  analogous to the proof of  theorem \ref{gen2}.
\qed

\section{Genus three fibrations on some surfaces with $p_g = 3, q = 0$}
The present section is meant to show, via a concrete example, how the 
structure theorem
for genus three fibrations \ref{strucgen3} can be applied.

Assume that $S$ is a minimal surface of general type with 
$p_g = 3, q = 0$,
and with $K_S^2 = 2 + d$. We  make several simplifying assumptions:

{\bf Assumption I :} we have a genus $3$ fibration $f : S \ra \PP^1$ 
( $B \cong \PP^1$ since $q=0$)
such that, for a general fibre $F$, the
restriction map  $ H^0(K_S) \ra H^0(\omega_F)$ is surjective.

It follows then that $V_1 \cong \oplus_1^3 \hol_{\PP^1}(2)$.
Our standard formulae read out as
$$K_S^2 = 2 + deg (\tau) \Rightarrow d = deg (\tau) .$$

{\bf Assumption II :} $d \leq 6$ and $\sigma_2$ yields a general 
extension class
\begin{equation}\label{sigma2pg=3}
0 \ra \oplus_1^6 \hol_{\PP^1}(4)
\ra \left( \oplus_1^d\hol_{\PP^1}(5) \right)\oplus\left( 
\oplus_1^{6-d}\hol_{\PP^1}(4)\right)
\ra \hol_{\tau}
\ra 0.
\end{equation}

It follows that $S^2\Lambda^2V_1 \cong \oplus_1^6 \hol_{\PP^1}(8)$, 
$\L'_4 \cong \hol_{\PP^1}(d+6)$ and
\begin{equation}\label{S2V2pg=3}
S^2(V_2)={6-d \choose 2} \hol_{\PP^1}(8) \oplus d(6-d)\hol_{\PP^1}(9) 
\oplus {d \choose 2}
\hol_{\PP^1}(10)
\end{equation}

\begin{prop}
Assume $d \leq 3$. Then $X$ is the complete intersection of $W$ with 
a relative quadric
$Z \in |\hol_{\PP(V_2)}(2) \otimes \pi^*(\hol_{\PP^1}(-6-d))|$, 
$\pi:\PP(V_2) \rightarrow \PP^1$ being
the natural projection.
\end{prop}

\Proof Tensoring the exact sequence
$$
0 \ra S^2(\Lambda^2(V_1))
\ra S^2(V_2)
\ra \tilde{V}_4
\ra 0.
$$
by ${\L'}_4^{-1}$ we see that the restriction map
$$H^0(S^2(V_2) \otimes {\L'}_4^{-1}) \rightarrow H^0(\tilde{V}_4 
\otimes {\L'}_4^{-1})$$
is surjective as soon as $H^1(S^2(\Lambda^2(V_1)) \otimes 
{\L'}_4^{-1})=0$. The result follows then from
$S^2(\Lambda^2(V_1)) \otimes {\L'}_4^{-1} \cong \oplus_1^6 \hol_{\PP^1}(2-d)$.
\qed
\begin{teo}
Let $d\leq 3$. Let $\sigma_2$ be a general homomorphism as in the 
exact sequence (\ref{sigma2pg=3}) and
let $\zeta$ be a general section of the vector bundle 
$S^2(V_2)\otimes {\L'}_4^{-1}$ (cf.
(\ref{S2V2pg=3})). Denote by $Z$ the divisor of $\zeta$ on $\PP(V_2)$ 
and set $X=W \cap Z$.
Then $X$ is a smooth minimal surface of general type with $p_g=3$, 
$q=0$ and $K^2_X=2+d$. We obtain in
this way a subvariety of the moduli space of dimension $33-2d$.
\end{teo}
\Proof
Let $U_1=\CC x_1 \oplus \CC x_2 \oplus  \CC x_3$ and $V_1 =U_1 
\otimes_{\CC} \hol_{\PP^1}(2)$.
Let $q_1, \ldots, q_6$ be a general basis of $Sym^2(U_1)$. Let 
$\sigma_2$ be given, in the
corresponding basis of $Sym^2(V_1)$, by the diagonal matrix whose 
entries are the first $6$ elements of a
sequence constituted by $1$ repeated $6-d$ times followed by $t_0$, 
$t_1$, $t_0-t_1$.

If $q_4,q_5,q_6$ yield irreducible conics, $W$ is a $\PP^2$ bundle 
except possibly over the points
$(0,1,\infty)$ where the fibre can be the cone over a rational normal 
quartic curve. The singular
points of $W$ are then exactly the cited vertices.

For $d\leq 2$, $S^2(V_2)\otimes {\L'}_4^{-1}$ is globally generated 
and therefore the linear system
cut by $Z$ on $W$ has no base points. So the smoothness assertion is
proved.

If $d=3$ we set $U_2:=\CC w_1 \oplus \CC w_2 \oplus  \CC w_3$, 
$U'_2:=\CC z_1 \oplus \CC z_2 \oplus  \CC
z_3$ and $V_2=(U_2 \otimes_{\CC} \hol_{\PP^1}(4)) \oplus (U'_2 
\otimes_{\CC} \hol_{\PP^1}(5))$. It
follows easily that in this case the base locus of the linear system 
of $Z$ is the $\PP^2$ bundle of
equations $\{z_1=z_2=z_3=0\}$.

Assume now that $q_4$, $q_5$, $q_6$ are smooth conics and that $q_4 
\cap q_5 \cap q_6 = \emptyset$.
Then we see that $W$ does not intersect the base locus of $Z$, and 
that $X$ is smooth.

The minimality of $X$ follows since $H^0(K_X)$ is given by 
$x_1,x_2,x_3$ and restricts to the canonical
system of each fibre. By the local equation that we have one sees 
that this last is base point free.

We compute the dimension of the  subvariety of the moduli 
space we have constructed. The choice of $\tau$ gives
$d$ parameters, the extension class $\xi$ (yielding $\sigma_2$) gives 
$6d-d=5d$ parameters since we
have to quotient by $Aut(\hol_{\tau})$, the linear system cut by $Z$ 
on $W$ has dimension $44-8d$. We
have therefore $44-2d$ parameters and we have to subtract $3$ for the 
automorphisms of $\PP^1$
and $8$ for those of $V_1$.
\qed
\begin{oss}
Minimal surfaces with $ K_S^2 = 2, p_g=3$ are well known since the 
time of Noether to be the double
covers of $\PP^2$ branched on a curve $D$ of degree $8$. Hence, their 
moduli space is
unirational of dimension $44 - 8 = 36$. The existence of 
such a genus three
fibration is equivalent to the existence of a one dimensional family
of quartic curves $C_t$
(images of the fibres) which are everywhere tangent to $D$. Thus 
 $D$ is a
square modulo each $C_t$, equivalently the equation of $D$ can be written as a
symmetrical determinant $ det \begin{pmatrix}A & B\\B & C \end{pmatrix} = 0 $. 
We obtain in this way only a 
$33$-dimensional
family, since $ 3 \times 15 - dim GL(2) - dim \PP GL(3) = 33$.

This shows how our above examples yield proper subvarieties of the 
moduli space.
\end{oss}

{\em Authors' addresses} \\
\noindent
  Prof. Dr. F. Catanese \\
Lehrstuhl Mathematik VIII \\
Universit\"at Bayreuth \\
Universit\"atsstr. 30, D-95447 Bayreuth\\

\noindent
Dr. R. Pignatelli \\
Dipartimento di Matematica \\
Universit\`a di Trento\\
Via Sommarive, 14, I- 38050 Povo (TN)
\end{document}